\documentclass[12pt]{article}  

\usepackage{amsfonts}
\usepackage{amssymb}

\newcommand{\cC}{{\cal C}}
\newcommand{\cD}{{\cal D}}
\newcommand{\cE}{{\cal E}}
\newcommand{\cF}{{\cal F}}
\newcommand{\cG}{{\cal G}}
\newcommand{\cH}{{\cal H}}
\newcommand{\cI}{{\cal I}}

\newcommand{\cL}{{\cal L}}
\newcommand{\cM}{{\cal M}}

\newcommand{\cO}{{\cal O}}
\newcommand{\cP}{{\cal P}}
\newcommand{\cR}{{\cal R}}

\newcommand{\gS}{{\Sigma}}
\newcommand{\gL}{{\Lambda}}
\newcommand{\gD}{{\Delta}}
\newcommand{\gO}{{\Omega}}
\newcommand{\gt}{{\theta}}
\newcommand{\gs}{{\sigma}}
\newcommand{\gG}{{\mit\Gamma}}

\newcommand{\Sym}{{\rm Sym}}
\newcommand{\rank}{{\rm rank}}
\newcommand{\op}{{\rm op}}
\newcommand{\coker}{{\rm coker}}
\newcommand{\ord}{{\rm ord}}

\renewcommand{\ker}{{\rm ker}}

\newcommand{\Aut}{{\rm Aut}}
\newcommand{\id}{{\rm id}}
\newcommand{\Spec}{{\rm Spec}}
\newcommand{\Proj}{{\rm Proj}}

\newcommand{\Ind}{{\rm Ind}}
\newcommand{\Res}{{\rm Res}}
\newcommand{\sgn}{{\rm sgn}}

\newcommand{\ZZ}{{\mathbb Z}}
\newcommand{\CC}{{\mathbb C}}
\newcommand{\FF}{{\mathbb F}}
\newcommand{\RR}{{\mathbb R}}
\newcommand{\NN}{{\mathbb N}}

\newcommand{\PP}{{\mathbb P}}
\newcommand{\QQ}{{\mathbb Q}}

\renewcommand{\ss}{\scriptstyle}
\newcommand{\ra}{\rightarrow}
\def\rightepi{{\longrightarrow \kern-0.7em \rightarrow}}
\newcommand{\oplusm}{\mathop{\oplus}\limits}
\newcommand{\otimesm}{\mathop{\otimes}\limits}
\newcommand{\HXl}{{\cH_{X,l}}}
\newcommand{\semidirect}{\rtimes}
\newcommand{\doublearrow}{\ra}
\newcommand{\notteilt}{{\,\not{\kern-0.075em|}\,}}

\begin{document}

\vspace*{15ex}

\begin{center}

{\LARGE\bf Riemann-Roch for Tensor Powers}\\
\bigskip
by\\
\bigskip
{\sc Bernhard K\"ock}

\end{center}

\bigskip

\section*{Introduction}

Let $f: X \ra Y$ be a projective morphism between noetherian schemes
$X$ and $Y$. Then, 
mapping a locally free $\cO_X$-module $\cF$ to the 
Euler characteristic $\sum_{i \ge 0}(-1)^i [R^if_*(\cF)]$ induces,
under certain additional assumptions, a push-forward homomorphism
\[f_*: K_0(X) \ra K_0(Y)\]
between the Grothendieck groups of $X$ and $Y$ (see \cite{FL}). The
fundamental task in Riemann-Roch theory is to compute this 
homomorphism $f_*$.

For any locally free $\cO_X$-module $\cF$, the symmetric group $\gS_l$
acts on the $l$-th tensor power $\cF^{\otimes l}$ by permuting the
factors. Using the ``binomial theorem''
\[(\cF \oplus \cG)^{\otimes l} \cong \oplusm_{i=0}^l 
\Ind_{\gS_i \times \gS_{l-i}}^{\gS_l} (\cF^{\otimes i} \otimes 
\cG^{\otimes (l-i)}),\]
the association $\cF \mapsto \cF^{\otimes l}$ can be canonically 
extended to a map
\[\tau^l: K_0(X) \ra K_0(\gS_l, X)\]
from $K_0(X)$ to the Grothendieck group $K_0(\gS_l, X)$ of all
locally free $\gS_l$-modules on $X$ (see section 1). We call $\tau^l$
the $l$-th tensor power operation. Similarly, we have an external
tensor power operation
\[\tau^l: K_0(X) \ra K_0(\gS_l, X^l)\]
where $X^l$ denotes the $l$-fold fibred product of $X$ with itself
over $Y$. 

In this paper, we study the following Riemann-Roch problem: How does
$\tau^l$ behave with respect to $f_*$? We give the following answers
to this problem.

{\bf Theorem A} (K\"unneth formula). Let $f$ be flat. Then the
following diagram commutes:
\[\begin{array}{ccc}
K_0(X) & \stackrel{\tau^l}{\longrightarrow } & K_0(\gS_l, X^l)\\
\\
{\ss f_*} \downarrow \phantom{\ss f_*} && 
\phantom{\ss f^l_*}\downarrow {\ss f^l_*} \\
\\
K_0(Y)&  \stackrel{\tau^l}{ \longrightarrow} & K_0(\gS_l, Y).
\end{array}\]
Here, $f^l_*$ denotes the equivariant push-forward homomorphism 
associated with the projection $f^l : X^l \ra Y$.

In the following theorems, let $\HXl$ denote the $\gS_l$-module
on $X$ of rank $l-1$ given by the short exact sequence
\[0 \ra \HXl \ra \cO_X[\{1,\ldots, l\}] \,\, \stackrel{\gS}{\ra}\,\,
\cO_X \ra 0\]
in which $\cO_X[\{1, \ldots, l\}]$ denotes the direct sum of $l$
copies of $\cO_X$ together with the obvious $\gS_l$-action and 
$\gS$ denotes the summation map.
Furthermore, we set $\lambda_{-1}(\cF) := \sum_{i \ge 0} (-1)^i 
[\gL^i(\cF)]$ for any locally free $\gS_l$-module $\cF$ on $X$.

{\bf Theorem B} (Riemann-Roch formula for closed immersions). Let $f$
be a regular closed immersion with conormal sheaf $\cC$. Then the
following diagram commutes:
\[\begin{array}{ccc}
K_0(X) & \stackrel{\lambda_{-1}(\cC \otimes \HXl) \cdot \tau^l}
{\longrightarrow} & K_0(\gS_l, X)\\
\\
{\ss f_*} \downarrow \phantom{\ss f_*} && 
\phantom{\ss f_*} \downarrow {\ss f_*} \\
\\
K_0(Y) & \stackrel{\tau^l} {\longrightarrow} & K_0(\gS_l, Y).
\end{array}\]

In the following two theorems, let $l$ be a prime and let $\tau^l$
denote also the composition
\[K_0(X) \,\, \stackrel{\tau^l}{\ra} \,\, K_0(\gS_l,X) \,\,
\stackrel{\Res^{\gS_l}_{C_l}}{\ra} \,\, K_0(C_l,X) \,\,
\stackrel{{\rm can}}{\rightepi} \,\, K_0(C_l,X)/(\cO_X[C_l])\]
where $C_l$ denotes the cyclic subgroup of $\gS_l$ generated by the
cycle $\langle 1, \ldots, l \rangle$.

{\bf Theorem C} (Riemann-Roch formula without denominators for 
principal $G$-bundles). Let $f$ be a principal $G$-bundle for some finite
group $G$ with $l \notteilt \ord(G)$. Then the following diagram commutes:
\[\begin{array}{ccc}
K_0(X) & \stackrel{\tau^l}{\longrightarrow} &
K_0(C_l,X)/(\cO_X[C_l])\\
\\
{\ss f_*} \downarrow \phantom{\ss f_*} && 
\phantom{\ss \bar{f}_*}\downarrow
{\ss \bar{f}_*} \\
\\
K_0(Y) & \stackrel{\tau^l}{\longrightarrow} & K_0(C_l,Y)/(\cO_Y[C_l]).
\end{array}\]

{\bf Theorem D} (Riemann-Roch formula for a general $f$). The
following diagram commutes:
\[\begin{array}{ccc}
K_0(X) & \stackrel{\lambda_{-1}(T_f^\vee \cdot \HXl)^{-1} \cdot
\tau^l}{\longrightarrow} & K_0(C_l,X)[l^{-1}]/(\cO_X[C_l]) \\
\\
{\ss f_*} \downarrow \phantom{\ss f_*} && 
\phantom{\ss \bar{f}_*} \downarrow {\ss \bar{f}_*} \\
\\
K_0(Y) & \stackrel{\tau^l}{\longrightarrow} & K_0(C_l,Y)[l^{-1}]/
(\cO_Y[C_l]).
\end{array}\]
Here, $T_f^\vee$ denotes the cotangential element associated with
$f$.

The proof of Theorem A (see section 2) essentially consists of a 
lengthy reduction to the ``base change isomorphism''
\[f^l_*(\cF_1 \boxtimes \ldots \boxtimes \cF_l) \cong 
f_*(\cF_1) \otimes \ldots \otimes f_*(\cF_l).\]
In this reduction, we will not use      
the K\"unneth spectral sequences developed in
\cite{EGA} but only some standard $K$-theoretical results of 
\cite{Q}. For the proof of Theorem B (see Theorem 4.1),
we may assume by using the deformation to the normal cone that 
$f$ is a so-called elementary embedding. In this case, the proof is
based on the multiplicativity of $\tau^l$ (see Proposition 1.7) and on
the formula 
\[\tau^l(\lambda_{-1}(\cF)) = \lambda_{-1}(\cF) \cdot 
\lambda_{-1}(\cF \otimes \HXl)\]
which is proved in Corollary 1.11. For the proof of Theorem C and 
Theorem~D (for smooth f), we generalize a new idea of Nori (see 
\cite{Ra}): In section 3, we prove the formula
\[\gD_*\left(\lambda_{-1}(\gO_{X/Y} \otimes \HXl)^{-1}\right) = 1 \quad
{\rm in} \quad K_0(C_l,X^l)[l^{-1}]/(\cO_{X^l}[C_l])\]
where $\gD: X \ra X^l$ denotes the diagonal. (This formula can also
be formulated for (non-projective) smooth morphisms $f$ between {\em
affine} schemes and can then even be strengthened to a $\gS_l$-version,
see Remark 3.8). Together with the K\"unneth formula, this formula
implies Theorem~D for smooth $f$ (see Theorem 4.2). 
(A $G$-equivariant version of)
Theorem C will be proved in Theorem 4.9. This proof is based on
the following ($G$-equivariant) version of the formula mentioned above:
\[\gD_*(1) = 1 \quad {\rm in} \quad K_0(C_l \times G, X^l)/ 
\Ind_{\{1\}}^{C_l}(K_0(G,X^l))\]
(see the proof of Theorem 4.9).

If $l$ is a prime, we have the following fundamental relation between
$\tau^l$ and the $l$-th Adams operation $\psi^l$: For all $x\in
K_0(X)$, we have:
\[\tau^l(x) = \psi^l(x) \quad {\rm in} \quad K_0(C_l, X)/(\cO_X[C_l])\]
(see Proposition 1.13). Furthermore, the element $\lambda_{-1}(\cF
\otimes \HXl)$ equals the $l$-th Bott element $\gt^l(\cF)$ in
$K_0(C_l, X)/(\cO_X[C_l])$ (see Proposition 3.2). Hence, Theorem~D
already follows from the classical Adams-Riemann-Roch formula for $f$
(see \cite{FL}). Conversely, Theorem~D implies the Adams-Riemann-Roch
formula if $X$ and $Y$ are $\CC$-schemes (see section 4). Thus, if
in addition $f$ is smooth, Nori's idea yields a new proof of the
Adams-Riemann-Roch-formula. If $f$ is a principal $G$-bundle such that
$l \notteilt \ord(G)$, Theorem C implies a version without denominators
of the (equivariant) Adams-Riemann-Roch formula 
$\psi^l \circ f_* = f_* \circ
\psi^l$ (see Corollary 4.10). One should be able to prove such
formulas for arbitrary \'etale morphism similarly to the paper 
\cite{FM} of Fulton and MacPherson.

Using Grayson's construction in \cite{GrEx}, one can
define (external) tensor power operations also for higher $K$-groups
(see section 1). Then, Theorem~A and Theorem C hold for higher
$K$-groups, too. If the conjecture on shuffle products in \cite{KoSh}
is true, then Theorem B and Theorem~D hold for higher $K$-groups, too.
Theorem~D further holds if $f$ is smooth or if the relation
``$\psi^l = \tau^l$'' holds also for higher $K$-groups.

In characteristic $p$, the $p$-th Adams operation is identical to
the pull-back homomorphism associated with the absolute Frobenius
morphism. This well-known relation can be considered as a substitute for
the relation ``$\psi^p = \tau^p$'' which is
extremely weak in characteristic $p$; 
for instance, we have $K_0(C_p, \FF_p,)/
(\FF_p[C_p])$ $\cong \ZZ/p\ZZ$, and Theorem~D for $Y=\Spec(\FF_p)$ and
$l=p$ even becomes trivial. Starting from this analogy, we develop
and investigate the following question in section 5, which may be
considered as the analogue of the formula in Theorem 3.1 mentioned
above: Let $f$ be smooth and let $F:X\ra X_Y$ denote the relative 
Frobenius morphism. Does then the following formula hold:
\[F_*(\gt^p(\gO_{X/Y})^{-1})=1 \quad {\rm in} \quad
K_0(X_Y)[p^{-1}] ?\]

{\bf Acknowledgments}. I would like to thank S.\ P.\ Dutta, R.\ H\"ubl, and 
G.\ Seibert for discussions about this question.
Furthermore, I would like to thank D.\ Grayson for his warm
hospitality during my stay at the University of Illinois at
Urbana-Champaign where the final part of this work was done.

\bigskip

{\bf Notations}. Let $G$ be an (abstract) group and $Y$ a
scheme. Then, a $G$-scheme over $Y$ is a $Y$-scheme $X$ together with
a homomorphism $G \ra \Aut_Y(X)^\op$ of groups; i.e., $G$ acts on
$X$ from the right by $Y$-automorphisms. A $G$-module on the
$G$-scheme $X$  is an $\cO_X$-module $\cF$ together with 
homomorphisms $g: g^*(\cF) \ra \cF$, $g\in G$, which satisfy the
usual homomorphism properties; i.e., $G$ acts on $\cF$ from the
left. The exact category of all locally free $\cO_X$-modules 
(respectively, of all locally free $G$-modules on $X$) of finite rank 
is denoted by $\cP(X)$ (respectively, by $\cP(G,X)$). The
corresponding higher $K$-groups (see \cite{Q}) are denoted by
$K_q(X)$ and $K_q(G,X)$, $q\ge 0$. Forgetting the $G$-structure
yields canonical homomorphisms $K_q(G,X) \ra K_q(X)$, $q\ge 0$. If
$G$ acts on $X$ trivially, we have natural maps $K_q(X) \ra
K_q(G,X)$, $q\ge 0$, which are right invers to the forgetful maps.\\
For any subgroup $H$ of $G$ and for any $H$-module $\cF$ on the 
$G$-scheme $X$, let $\Ind_H^G(\cF)$ denote the $G$-module 
$\oplus_{r\in R}\, r^*(\cF)$ on $X$ where $R \subseteq G$ is a 
system of representatives for $G/H$ and where $G$ acts on 
$\Ind_H^G(\cF)$ as follows: For any pair of elements $g \in G$, 
$r \in R$, let $s\in R$ and $h\in H$ be the uniquely determined 
elements such that $gr = sh$. Then $g$ acts on the direct summand
$r^*(\cF)$ of $\Ind_H^G(\cF)$ by virtue of the composition
\[g^*(r^*(\cF)) = (gr)^*(\cF) = (sh)^*(\cF) = s^*(h^*(\cF))
\,\, \stackrel{s^*(h)}{\longrightarrow} \,\, s^*(\cF).\]
It is easy to prove that $\Ind_H^G(\cF)$ does not depend on the 
chosen system $R$ of representatives (up to a canonical
isomorphism).\\
For any $l \ge 0$, $\gS_l$ denotes the group of permutations of
the set $I_l := \{1, \ldots, l\}$. For any $Y$-scheme $X$, the
$l$-fold fibred product $X^l := X \times_Y \ldots \times_Y X$ is
considered as a $\gS_l$-scheme by
\[\sigma(x_1, \ldots, x_l) := (x_{\sigma(1)}, \ldots, x_{\sigma(l)})\]
(for $x_1, \ldots, x_l \in X$ and $\gs \in \gS_l$). For any 
$\gS_l$-scheme $X$, $\cO_X[I_l]$ denotes the $\cO_X$-module 
$\oplusm^l \cO_X$ together with the obvious $\gS_l$-action. The
$\gS_l$-module $\HXl$ on any $\gS_l$-scheme $X$ is defined by
the short exact sequence
\[0 \ra \HXl \ra \cO_X[I_l] \,\, \stackrel{\gS}{\ra} \,\, 
\cO_X \ra 0.\]
The cyclic subgroup of $\gS_l$ generated by the cycle $c :=
\langle 1, \ldots, l \rangle \in \gS_l$ is denoted by $C_l$.
$\HXl$ is isomorphic to $\coker(\gD: \cO_X \,\, \ra
\,\, \cO_X[I_l])$ as a $C_l$-module and, if $l$ is invertible on
$X$, even as a $\gS_l$-module.

\bigskip

\section*{\S 1 (External) Tensor Power Operations}

Let $f: X \ra Y$ be a morphism between noetherian schemes. In this
section, we will construct natural external tensor power operations
\[\tau^l: K_q(X) \ra K_q(\gS_l, X^l), \quad q \ge 0, \quad l \ge 0,\]
such that, for $q =0$, the class $[\cF]$ of any locally free 
$\cO_X$-module $\cF$ is mapped to the class $[\cF^{\boxtimes l}]$ of
the external tensor power $\cF^{\boxtimes l}$ where $\gS_l$ acts on
$\cF^{\boxtimes l}$ by permuting the factors. In particular, for $f= \id$,
we obtain tensor power operations
\[\tau^l: K_q(Y) \ra K_q(\gS_l, Y), \quad  q\ge 0, \quad l \ge 0. \]
We will show in Proposition 1.7 that $\tau^l$ is multiplicative (on
Grothendieck groups). Furthermore, we will prove that, for any locally
free $\cO_Y$-module $\cF$ the element $\tau^l(\lambda_{-1}([\cF]))$
is divisible by the element $\lambda_{-1}([\cF])$ in $K_0(\gS_l, Y)$
(see Corollary 1.11). If $l$ is a prime, we finally establish a
fundamental relation between the tensor power operation
\[K_q(Y) \,\, \stackrel{\tau^l}{\longrightarrow} \,\,
K_q(\gS_l, Y) \,\, \stackrel{\Res_{C_l}^{\gS_l}}{\longrightarrow}
\,\, K_q(C_l, Y)\]
and the Adams operation $\psi^l$ on $K_q(Y)$ (see Proposition 1.13)
and we will explicitly compute the operation $\tau^l: K_1(\CC) \ra
K_1(C_l, \CC)$ (see Proposition 1.14).

We now begin to construct $\tau^l$. For any $l\ge 0$, let 
$F_l(\cP(X))$ denote the category of all chains $\cF_1 
\hookrightarrow \ldots \hookrightarrow \cF_l$ of admissible 
monomorphisms in $\cP(X)$. Here, a monomorphism $\cF \hookrightarrow
\cG$ in $\cP(X)$ is called admissible if the quotient $\cG/\cF$ is
locally free. We have functors
\[\times: \cP(\gS_k, X^k) \times \cP(\gS_l, X^l) \ra 
\cP(\gS_{k+l}, X^{k+l}), \quad (\cF, \cG) \mapsto 
\Ind_{\gS_k \times \gS_l}^{\gS_{k+l}}(\cF \boxtimes \cG),\]
(for all $k,l \ge 0$) and
\[ F_l(\cP(X)) \ra \cP(\gS_l, X^l), \quad (\cF_1 \hookrightarrow
\ldots \hookrightarrow \cF_l) \mapsto \sum_{\gs \in \gS_l} 
\cF_{\gs(1)} \boxtimes \ldots \boxtimes \cF_{\gs(l)},\]
(for all $l \ge 0$). Here, $\cF \boxtimes \cG$ denotes the external
tensor product $p^*(\cF) \otimes q^*(\cG)$ for $p: X^{k+l} \ra X^k$ and
$q: X^{k+l} \ra X^l$ the projections to the first $k$ and last
$l$ components, respectively. It is considered as a $\gS_k \times
\gS_l$-module on the $\gS_{k+l}$-scheme $X^{k+l}$ in the obvious way.
Similarly, the sum 
$\sum_{\gs \in \gS_l} \cF_{\gs(1)} \boxtimes \ldots
\boxtimes \cF_{\gs(l)}$ 
is considered as a $\gS_l$-submodule of 
the external tensor power $\cF_l^{\boxtimes l}$. It is easy, but tedious,
to prove that this sum of submodules is a locally free $\cO_X$-module
and that these functors satisfy the axioms (E1) through (E5) in section 7
of Grayson's paper \cite{GrEx} (see also section 2 in \cite{KoSh},
last example in \cite{GrEx}, and Remark 1.6).

{\bf Remark 1.1}. Similarly to the last example in \cite{GrEx}, one
can define natural power operations $F_l(\cP(\gS_k, X^k)) \ra 
\cP(\gS_{kl}, X^{kl})$ for all $k \ge 1$. However, these operations
can essentially be reduced to the case $k = 1$ (see Remark 1.12).
So, we will consider only the case $k=1$ in order to 
avoid complicated notations. 

Grayson's construction in section 7 of \cite{GrEx} yields maps
\[\tau^l: K_q(X) \ra K_q(\gS_l, X^l), \quad q \ge 0, \quad l\ge 0,\]
which we will call {\em external tensor power operations}. 
If $f= \id$, we in particular obtain maps
\[\tau^l: K_q(Y) \ra K_q(\gS_l, Y),\quad  q\ge 0, \quad l \ge 0,\]
which we will call {\em tensor power operations}. In the following
lemma, we describe $\tau^l$ on Grothendieck groups.

{\bf Lemma 1.2}. For any $l\ge 0$ and any locally free $\cO_X$-modules
$\cE$, $\cF$, we have in $K_0(\gS_l, X^l)$:
\[\tau^l([\cE]-[\cF]) = \sum_{{a \ge 0, b_1, \ldots, b_u \ge 1} \atop
{a+b_1+\ldots + b_u = l}} (-1)^u\left[\Ind_{\gS_a {\times} \gS_{b_1} {\times}
\ldots {\times} \gS_{b_u}}^{\gS_l} (\cE^{\boxtimes a} {\boxtimes} 
\cF^{\boxtimes b_1} {\boxtimes} \ldots {\boxtimes} \cF^{\boxtimes b_u})\right].\]
In particular, we have
\[\tau^l([\cE]) = [\cE^{\boxtimes l}] \quad {\rm in} \quad
K_0(\gS_l, X^l).\]

{\bf Proof}. This follows from section 8 in \cite{GrEx}.

For any $i, j \ge 0$, the functor $\times $ obviously induces a 
bilinear product
\[\times: K_0(\gS_i, X^i) \times K_0(\gS_j, X^j) \ra 
K_0(\gS_{i+j}, X^{i+j}) \]
which we call {\em cross product}. The abelian group 
$\prod_{l\ge 0} K_0(\gS_l, X^l)$ together with the cross product
\[\begin{array}{cccl}
\times: & \prod_{i\ge 0} K_0(\gS_i, X^i) \times 
\prod_{j\ge 0} K_0(\gS_j, X^j) & \ra & \prod_{l\ge 0} 
K_0(\gS_l, X^l) \\
& ((x_i)_{i\ge 0}, (y_j)_{j\ge 0}) & \mapsto & 
(\sum_{i+j = l} x_i \times y_j)_{l\ge 0}
\end{array}\]
is then a commutative ring with $1$ (see also Proposition 2.1 on 
p.\ 12 in \cite{Ho}).

{\bf Corollary 1.3}. The map
\[\tau := (\tau_l)_{l\ge 0} : (K_0(X), +) \ra 
(\prod_{l\ge 0} K_0(\gS_l, X^l), \times)\]
is a homomorphism.

{\bf Proof}. For any locally free $\cO_X$-modules $\cE$, $\cF$, we 
have
\[\tau([\cE] - [\cF]) = ([\cE^{\boxtimes l}])_{l \ge 0} \times
([\cF^{\boxtimes l}])_{l \ge 0}^{-1} = \tau([\cE]) \times
\tau([\cF])^{-1}\]
by section 2 in \cite{GrEx} and by Lemma 1.2. Furthermore, for 
all $l \ge 0$, we have a canonical isomorphism 
\[(\cE \oplus \cF)^{\boxtimes l} \cong \oplusm_{i= 0}^l 
\Ind_{\gS_i \times \gS_{l-i}}^{\gS_l}(\cE^{\boxtimes i} \boxtimes 
\cF^{\boxtimes (l-i)})\]
of $\gS_l$-modules on $X^l$; i.e., we have 
\[\tau([\cE \oplus \cF]) = \tau([\cE]) \times \tau([\cF]).\]
This immediately proves Corollary 1.3.

{\bf Proposition 1.4}. For any $l\ge 0$ and $1\le i \le l-1$, the
following diagram commutes:
\[\begin{array}{ccc}
K_0(X) & \stackrel{\tau^l}{\longrightarrow} & K_0(\gS_l, X^l)\\
\\
{\ss (\tau^i, \tau^{l-i})} \downarrow \phantom{\ss (\tau^i, \tau^{l-i})} &&
\phantom{\ss \Res_{\gS_i \times \gS_{l-e}}^{\gS_l}} \downarrow 
{\ss \Res_{\gS_i \times \gS_{l-i }}^{\gS_l}}\\
\\
K_0(\gS_i, X^i) \times K_0(\gS_{l-i}, X^{l-i}) & 
\stackrel{\boxtimes}{\longrightarrow} & 
K_0(\gS_i \times \gS_{l-i}, X^i \times X^{l-i}).
\end{array}\]
If Conjecture' on p.\ 289 in \cite{KoSh} is true, then the
corresponding diagrams of higher $K$-groups commute, too.

{\bf Proof.} First, we prove the assertion for Grothendieck
groups. For this, we define a double cross product $\times \times$ on
$\prod_{i,j} K_0(\gS_i \times \gS_j, X^{i+j})$ by
\[[\cE] {\times} {\times}  [\cF] :=
\left[ 
\Ind_{\gS_{i_1}\times \gS_{j_1} \times \gS_{i_2} \times \gS_{j_2}}
^{\gS_{i_1+i_2}\times \gS_{j_1+j_2}}(\cE \boxtimes \cF)\right]\]
for $\cE \in \cP(\gS_{i_1}{\times} \gS_{j_1}, X^{i_1 + j_1})$ and
$\cF \in \cP(\gS_{i_2} {\times} \gS_{j_2}, X^{i_2+j_2})$.
As in Proposition 2.1 on p.\ 12 in \cite{Ho}, one shows that the 
abelian group $\prod_{i, j} K_0(\gS_i\times \gS_j, X^{i+j})$
together with $\times \times$ is a commutative ring with $1$. Using
Mackey's subgroup theorem, one shows as in Theorem 1.2 on p.\ 8 in
\cite{Ho} that the restriction map
\[\Res: (\prod_{l\ge 0} K_0(\gS_l, X^l), \times) \ra 
(\prod_{i,j} K_0(\gS_i \times \gS_j, X^{i+j}), \times \times)\]
is a homomorphism. It follows from this, from Corollary 1.3, and from
the definition of $\times \times$ that all maps in the diagram 
\[\begin{array}{ccc}
K_0(X) & \stackrel{\tau}{\ra} & \prod_l K_0(\gS_l, X^l)\\
\\
{\ss (\tau, \tau)} \downarrow \phantom{\ss (\tau, \tau)} &&
\phantom{\ss \Res} \downarrow {\ss \Res}\\
\\
\prod_i K_0(\gS_i, X^i) \times \prod_j K_0(\gS_j, X^j) & 
\stackrel{\boxtimes}{\ra} & \prod_{i,j} K_0(\gS_i \times \gS_j,
X^{i+j})
\end{array}\]
are homomorphisms. Thus, it suffices to show that, for all
locally free $\cO_X$-modules $\cE$, the following equality holds:
\[\Res( \tau ([\cE])) = \tau([\cE]) \boxtimes \tau([\cE])
\quad {\rm in} \quad \prod_{i,j} K_0(\gS_i \times \gS_j, X^{i+j}).\]
This equality follows from the obvious fact that,
for all $i, j$, the 
$\gS_i \times \gS_j$-modules $\cE^{\boxtimes i} \boxtimes 
\cE^{\boxtimes j}$ and $\Res_{\gS_i \times \gS_j}^{\gS_{i+j}}
(\cE^{\boxtimes (i+j)})$ on 
$X^{i+j}$ are isomorphic.\\
For the proof of the corresponding assertion for the higher
$K$-groups, we use the notations introduced in
\cite{KoSh}. Furthermore, we set $\cP_l := \cP(\gS_l, X^l)$ for all
$l \ge 0$. By construction (see \cite{GrEx} or \cite{KoSh}), the map
$\tau^l: K_q(X) \ra K_q(\gS_l, X^l)$ is induced by the 
continuous map
\[|G\cP_1| \, \, \stackrel{|\tau^l|}{\ra} \,\, 
|\cH^l\cP_l| \,\, \stackrel{\Xi^l}{\ra} \,\, |G^l\cP_l|.\]
Here, the first map is the geometric realization of a certain 
simplicial map $G\cP_1 \ra \cH^l\cP_l$ which is induced by the
above functors and which we denote by $\tau^l$, too. The second map
is defined in sections 5 and 7 in \cite{GrEx}. One easily sees
that the diagram
\[\begin{array}{ccc}
G\cP_1 & \stackrel{\tau^l}{\longrightarrow} & \cH^l\cP_l\\
\\
{\ss {}^{\boxtimes l}} \downarrow \phantom{\ss {}^{\boxtimes l}} &&
\phantom{\ss \Res_{\gS_1}^{\gS_l}} \downarrow 
{\ss \Res_{\gS_1}^{\gS_l}} \\
\\
G^l\cP(X^l) & \stackrel{\Sym}{\longrightarrow} &
\cH^l\cP(X^l)
\end{array}\]
of simplicial maps commutes (up to homotopy). Here, the simplicial 
map ${}^{\boxtimes l}$ is induced by the external tensor product,
and $\Sym := \Sym^{1,\ldots, 1}$ is the shuffle operation introduced in
(the last remark of) section 3 in \cite{KoSh}. Since the external 
tensor product and
the shuffle operation are associative in the obvious sense, the
diagram
\[\begin{array}{ccccc}
G\cP_1\\
\\
{\scriptstyle {({}^{\boxtimes i}, {}^{\boxtimes (l-i)})}} \downarrow
\phantom{\scriptstyle ({}^{\boxtimes i}, {}^{\boxtimes (l-i)})} 
& \phantom{\scriptstyle {}^{\boxtimes l}} 
\searrow {\scriptstyle {}^{\boxtimes l}}\\
\\
G^i\cP(X^i) \times G^{l-i}\cP(X^{l-i}) & 
\stackrel{\boxtimes}{\ra} & G^l \cP(X^l)\\
\\
{\scriptstyle \Sym\times \Sym} \downarrow
\phantom{\scriptstyle \Sym\times \Sym} 
&& & \phantom{\scriptstyle \Sym}\searrow 
{\scriptstyle \Sym}\\
\\
\cH^i\cP(X^i) \times \cH^{l-i}\cP(X^{l-i}) & 
\stackrel{\boxtimes}{\ra} & \cH^{i, l-i}\cP(X^l) &
\stackrel{\Sym^{i, l-i}}{\ra} & \cH^l\cP(X^l)
\end{array}\]
of simplicial maps commutes (up to homotopy). This implies that
the upper quadrangle in the diagram
\[\begin{array}{ccccccc}
|G\cP_1| &&& \stackrel{|\tau^l|}{\ra} &&& |\cH^l\cP_l|\\
\\
{\ss (|\tau^i|, |\tau^{l-i}|)} \downarrow 
\phantom{\ss (|\tau^i|, |\tau^{l-i}|)} &&&&& 
\phantom{\ss {\rm Res}} \swarrow {\ss {\rm Res}}
& \phantom{\ss \Xi^l} \downarrow {\ss \Xi^l} \\
\\
|\cH^i\cP_i| \times |\cH^{l-i} \cP_{l-i}| & 
\stackrel{\boxtimes}{\ra} & |\cH^{i,l-i}\cP_{i,l-i}| & 
\stackrel{|\Sym^{i,l-i}|}{\ra} & |\cH^l\cP_{i,l-i}| && |G^l\cP_l| \\
\\
{\ss \Xi^i \times \Xi^{l-i}} \downarrow 
\phantom{\ss \Xi^i \times \Xi^{l-i}} &&&
{\ss \Xi^{i, l-i}} \searrow \phantom{\ss \Xi^{i, l-i}}
& \phantom{\ss \Xi^l} \downarrow {\ss \Xi^l} &
\phantom{\ss {\rm Res}} \swarrow {\ss {\rm Res}}\\
\\
|G^i\cP_i| \times |G^{l-i}\cP_{l-i}| && 
\stackrel{\boxtimes}{\ra} && |G^l\cP_{i,l-i}|
\end{array}\]
commutes up to homotopy ($\cP_{i, l-i} := \cP(\gS_i \times
\gS_{l-i}, X^l)$). The lower left and the right quadrangle commute
for trivial reasons. Conjecture' on p.\ 289 in \cite{KoSh} implies
that the lower middle triangle commutes up to homotopy. Thus, the
exterior pentagon commutes, too. This proves Proposition 1.4 for 
higher $K$-groups.

{\bf Remark 1.5}. The conjecture on shuffle products mentioned above
has been proved in the case $i=1$ by Nenashev (see \cite{Ne}). Thus, 
Proposition 1.4 holds in the case $i=1$ also for higher $K$-groups. 
Furthermore, it follows from this (or already from the corollary on
p.\ 293 in \cite{KoSh}) that the diagram
\[\begin{array}{ccccc}
|G\cP_1| & \stackrel{|\tau^l|}{\ra} & |\cH^l\cP_l| & 
\stackrel{\Xi^l}{\ra} & |G^l\cP_l| \\
\\
& {\ss {}^{\boxtimes l}} \searrow \phantom{\ss {}^{\boxtimes l}}
&&\phantom {\ss \Res^{\gS_l}_{\gS_1}} \swarrow {\ss \Res^{\gS_l}_{\gS_1}}\\
\\
&& |G^l\cP(X^l)|
\end{array}\]
of continuous maps commutes up to homotopy; i.e., $\Xi^l \circ
|\tau^l|$ is an equivariant lift of the usual external tensor power 
operation ${}^{\boxtimes l}$. In particular, the composition
\[K_q(X) \,\, \stackrel{\tau^l}{\longrightarrow} \,\, 
K_q(\gS_l, X^l) \,\, \stackrel{\Res^{\gS_l}_{\gS_1}}{\longrightarrow}
\,\, K_q(X^l)\]
is identical to the map $x \mapsto p_1^*(x) \cdot \ldots \cdot
p_l^*(x)$ which, for $q\ge 1$, is the zero map since products on
higher $K$-groups are trivial.

{\bf Remark 1.6}. Without using the above functors, the
simplicial maps
\[\tau^l : G\cP(X) \ra \cH^l\cP(\gS_l, X^l), \quad l\ge 0,\]
(cf.\ proof of Proposition 1.4) and then the (external) tensor
power operations on higher $K$-groups can be constructed as follows:
Similarly to section 3 in \cite{KoAdPro}, one shows that the 
composition
\[G\cP(X) \,\, \stackrel{{}^{\boxtimes l}}{\longrightarrow} \,\,
G^l\cP(X^l) \,\, \stackrel{\Sym^{1, \ldots, 1}}{\longrightarrow}
\,\, \cH^l \cP(X^l)\]
can be canonically lifted to a simplicial map $\tau^l: 
G\cP(X) \ra \cH^l\cP(\gS_l, X^l)$. This construction of $\tau^l$ 
has the advantage that Grayson's axioms for power operations 
need not to be checked for the above functors.

{\bf Proposition 1.7}. \\
(a) For any $l\ge 0$ and any $x,y \in K_0(X)$, we have:
\[\tau^l(x \cdot y) = \tau^l(x) \cdot \tau^l(y) \quad {\rm in} 
\quad K_0(\gS_l, X^l).\]
(b) Let $q \ge 1$. For any $x \in K_0(X)$ of the form $x=[\cE]$
($\cE$ a locally free $\cO_X$-module) and for any $y \in K_q(X)$,
we have:
\[\tau^l(x \cdot y) = \tau^l(x) \cdot \tau^l(y) \quad {\rm in}
\quad K_q(\gS_l, X^l).\]
If Conjecture' on p.\ 289 in \cite{KoSh} is true, then this equality
holds for all $x\in K_0(X)$. 

{\bf Proof}. \\
(a) Using Corollary 1.3, Frobenius reciprocity, and
Proposition 1.4 successively, we obtain for all $x_1, x_2, y \in K_0(X)$ and
$l\ge 0$: 
\begin{eqnarray*}
\lefteqn{\tau^l(x_1 + x_2) \cdot \tau^l(y) =}\\
&& = \sum_{i=0}^l \Ind_{\gS_i \times \gS_{l-i}}^{\gS_l}
(\tau^i(x_1) \boxtimes \tau^{l-i}(x_2)) \otimes \tau^l(y) \\
&& = \sum_{i=0}^l \Ind_{\gS_i \times \gS_{l-i}}^{\gS_l}
\left((\tau^i(x_1) \boxtimes \tau^{l-i}(x_2)) \otimes 
\Res_{\gS_i \times \gS_{l-i}}^{\gS_l}(\tau^l(y))\right) \\
&& = \sum_{i=0}^l \Ind_{\gS_i \times \gS_{l-i}}^{\gS_l}
\left((\tau^i(x_1) \otimes \tau^i(y)) \boxtimes (\tau^{l-i}(x_2) \otimes
\tau^{l-i}(y))\right) \\
&& = l\mbox{-th component of } (\tau(x_1) \cdot \tau(y)) \times
(\tau(x_2) \cdot \tau(y)).
\end{eqnarray*}
This computation shows that the map
\[\begin{array}{ccl}
K_0(X) \times K_0(X) &\ra &(\prod_{l \ge 0} K_0(\gS_l, X^l), \times)\\
(x,y) &\mapsto &\tau(x) \cdot \tau(y) := 
(\tau^l(x) \cdot \tau^l(y))_{l\ge 0}
\end{array}\]
is bilinear. By Corollary 1.3, the map
\[K_0(X) \times K_0(X) \ra (\prod_{l\ge 0} K_0(\gS_l, X^l), \times),
\quad (x,y) \mapsto \tau(xy), \]
is bilinear, too. Hence, it suffices to show the equality
\[\tau([\cE]) \cdot \tau([\cF]) = \tau([\cE \otimes \cF]) \quad 
{\rm in} \quad \prod_{l\ge 0} K_0(\gS_l, X^l)\]
for all locally free $\cO_X$-modules $\cE$ and $\cF$. This equality 
immediately follows from the obvious fact that, for all $l \ge 0$,
the $\gS_l$-modules $\cE^{\boxtimes l} \otimes \cF^{\boxtimes l}$ 
and $(\cE \otimes \cF)^{\boxtimes l}$ are isomorphic.\\
(b) For any locally free $\cO_X$-modules $\cF_1, \ldots, \cF_u$ 
and for any $b_1, \ldots, b_u \in \NN$ with $b_1 + \ldots + b_u =l$,
we have a canonical isomorphism
\[\cE^{\boxtimes l} {\otimes} \Ind_{\gS_{b_1} {\times} \ldots {\times} 
\gS_{b_u}}^{\gS_l} (\cF_1^{\boxtimes b_1} {\boxtimes} \ldots {\boxtimes}
\cF_u^{\boxtimes b_u}) \cong
\Ind_{\gS_{b_1} {\times} \ldots {\times} \gS_{b_u}}^{\gS_l}
((\cE{\otimes} \cF_1)^{\boxtimes b_1} {\boxtimes} \ldots {\boxtimes}
(\cE {\otimes} \cF_u)^{\boxtimes b_u})\]
of $\gS_l$-modules on $X^l$ by Frobenius reciprocity. Similarly to
Proposition 7.2 on p.\ 306 in \cite{KoSh}, one deduces the first 
assertion in (b) from this. The second assertion then follows 
as in (a).

{\bf Remark 1.8}. One can find constructions and statements which
are similar to those made so far for Grothendieck groups already in
Hoffman's book \cite{Ho} (in particular, see Proposition 3.7
on p.\ 22 in \cite{Ho}). However, he uses the simplifying facts
that short exact sequences of vector bundles split and that, for
instance, $K_0(\gS_i, X^i) \otimes K_0(\gS_j, X^j)$ is isomorphic
to $K_0(\gS_i \times \gS_j, X^{i+j})$. Since such facts are not
available in our situation, Hoffman's proofs have been adapted 
to our situation and included into this paper for the reader's
convenience.

{\bf Proposition 1.9}. For any $l\ge 0$ and any locally free
$\cO_X$-module $\cF$, we have:
\[\tau^l(-[\cF]) = (-1)^l [\cF^{\boxtimes l}_\sgn] \quad
{\rm in} \quad K_0(\gS_l, X^l);\]
here, $\cF^{\boxtimes l}_\sgn := \cF^{\boxtimes l} \otimes 
\cO_{X^l, \sgn}$ denotes the tensor product of the $\gS_l$-module
$\cF^{\boxtimes l}$ on $X^l$ with the sign representation 
$\cO_{X^l,\sgn}$ of $\gS_l$. 

{\bf Proof}. First, let $l \ge 0$ be fixed. The element 
$N:= \sum_{i=1}^l [i]$ of the representation $\cO_{X^l}[I_l]$
is obviously $\gS_l$-invariant. Thus, the homomorphisms
\[d_i: \gL^i(\cO_{X^l}[I_l]) \ra \gL^{i+1} (\cO_{X^l}[I_l]),
\quad x \mapsto N \wedge x,\]
$i=0, \ldots, l-1$, are $\gS_l$-homomorphisms and define a complex
\[0 \ra \cO_{X^l} \,\, \stackrel{d_0}{\ra} \,\, \cO_{X^l}[I_l] 
\,\, \stackrel{d_1}{\ra} \,\, \ldots \,\, 
\stackrel{d_{l-2}}{\ra} \,\, \gL^{l-1}(\cO_{X^l}[I_l]) \,\,
\stackrel{d_{l-1}}{\ra} \,\, \gL^l(\cO_{X^l}[I_l]) \ra 0\]
of $\gS_l$-modules on $X^l$. Let $p_1: \cO_{X^l}[I_l] \ra 
\cO_{X^l}$ denote the projection onto the $[1]$-component. Then
one easily checks that the maps
\[\begin{array}{ccl}
\gL^{i+1}(\cO_{X^l}[I_l])& \ra &\gL^i(\cO_{X^l}[I_l])\\
f_1 \wedge \ldots \wedge f_{i+1} &\mapsto&
\sum_{a=1}^{i+1} (-1)^{a-1} p_1(f_a)\, f_1 \wedge \ldots \wedge
\hat{f}_a \wedge \ldots \wedge f_{i+1}
\end{array}\]
$i=1, \ldots, l$, are well-defined and form a (non-equivariant)
homotopy between the identity and the zero map on this complex. 
Hence, the complex is exact. For any $i=0, \ldots, l$, we have
an isomorphism
\[\gL^i(\cO_{X^l}[I_l]) \,\, \tilde{\ra} \,\, 
\Ind_{\gS_i \times \gS_{l-i}}^{\gS_l} (\cO_{X^i, \sgn} \boxtimes
\cO_{X^{l-i}})\]
of $\gS_l$-modules on $X^l$ which is defined as follows: For 
any $\rho(1) < \ldots < \rho(i)$ in $I_l$, the basis element
$[\rho(1)] \wedge \ldots \wedge [\rho(i)]$ of $\gL^i(\cO_{X^l}[I_l])$
is mapped to the basis element $[\rho]$ of $\Ind_{\gS_i \times 
\gS_{l-i}}^{\gS_l}(\cO_{X^i,\sgn} \boxtimes \cO_{X^{l-i}})$ where
$\rho$ denotes the $(i,l-i)$-shuffle permutation corresponding to
$\rho(1) < \ldots < \rho(i)$. Tensoring the above complex with
$\cF^{\boxtimes l}_\sgn$ and using Frobenius reciprocity, we obtain
the following exact sequence of $\gS_l$-modules on $X^l$:
\[0 \ra \cF^{\boxtimes l}_\sgn \ra 
\Ind_{\gS_1 \times \gS_{l-1}}^{\gS_l} (\cF \boxtimes 
\cF^{\boxtimes (l-1)}_\sgn) \ra \ldots \qquad \qquad \qquad\]
\[\qquad \qquad \qquad \ra 
\Ind_{\gS_{l-1} \times \gS_1}^{\gS_l} (\cF^{\boxtimes (l-1)}
\boxtimes \cF) \ra \cF^{\boxtimes l} \ra 0.\]
This shows that the element 
$((-1)^l [\cF^{\boxtimes l}_\sgn ])_{l\ge 0} \in \prod_{l\ge 0}
K_0(\gS_l, X^l)$ is an inverse of the element $\tau([\cF]) =
([\cF^{\boxtimes l}])_{l\ge 0}$ with respect to the cross product. Now,
Proposition 1.9 follows from Corollary 1.3.

As seen in Corollary 1.3, the construction of the external tensor
power operation $\tau^l$  is based on the equivariant version
\[(\cE \oplus \cF)^{\boxtimes l} \cong \oplusm_{i=0}^l 
\Ind_{\gS_i \times \gS_{l-i}}^{\gS_l} (\cE^{\boxtimes i} \boxtimes
\cF^{\boxtimes (l-i)})\]
of the binomial theorem $(a+b)^l = \sum_{i=0}^l {l \choose i} 
a^i b^{l-i}$. The following corollary contains the equivariant version
of the binomial theorem $(a-b)^l = \sum_{i=0}^l (-1)^{l-i}
{l \choose i} a^i b^{l-i}$.

{\bf Corollary 1.10}. For any $l\ge 0$ and any locally free 
$\cO_X$-modules $\cE$, $\cF$, we have:
\[\tau^l([\cE]-[\cF]) = \sum_{i=0}^l (-1)^{l-i}\left[
\Ind_{\gS_i \times \gS_{l-i}}^{\gS_l} (\cE^{\boxtimes i} \boxtimes
\cF^{\boxtimes (l-i)}_\sgn)\right] \quad {\rm in} \quad 
K_0(\gS_l, X^l).\]

{\bf Proof}. This follows from Corollary 1.3 and Proposition 1.9.

The next corollary will be an essential ingredient in the proof
of the Riemann-Roch formula for closed immersions (see Theorem 4.1).
Let $\lambda_{-1}(x)$ denote the element 
\[\lambda_{-1}(x) := \sum_{i\ge 0} (-1)^i \lambda^i(x) \in K\]
for any $x$ of finite $\lambda$-degree in a $\lambda$-ring $K$. 

{\bf Corollary 1.11}. Let $f=\id$. For any $l\ge 0$ and any locally
free $\cO_Y$-module $\cF$, we have:
\[\tau^l(\lambda_{-1}([\cF])) = \lambda_{-1}([\cF \otimes 
\cO_Y[I_l]]) = \lambda_{-1}([\cF \otimes \cH_{Y,l}]) \cdot
\lambda_{-1}([\cF]) \quad {\rm in} \quad K_0(\gS_l, Y).\]

{\bf Proof}. The second equality is clear since $\lambda_{-1}$ is
multiplicative. We first prove the first equation in the case
$\rank(\cF) = 1$. Then we have by Corollary 1.10:
\begin{eqnarray*}
\lefteqn{\tau^l(\lambda_{-1}([\cF])) = \tau^l([\cO_Y] - [\cF]) }\\
&& = \sum_{i=0}^l (-1)^{l-i} 
\left[\Ind_{\gS_i \times \gS_{l-i}}^{\gS_l}(\cF^{\otimes (l-i)}_\sgn)
\right] \\
&& = \sum_{i=0}^l (-1)^{l-i}\left[\cF^{\otimes (l-i)} \otimes 
\gL^{l-i}(\cO_Y[I_l])\right] \\
&& = \lambda_{-1}([\cF \otimes \cO_Y[I_l]]) \quad {\rm in}
\quad K_0(\gS_l,Y).
\end{eqnarray*}
(In the third equality, we have used the isomorphisms
\[\Ind_{\gS_i \times \gS_{l-i}}^{\gS_l} (\cO_Y \otimes \cO_{Y,\sgn})
\cong \gL^{l-i}(\cO_Y[I_l]), \quad i=0, \ldots, l,\]
established in the proof of Proposition 1.9.) Now let $\rank(\cF)$ be
arbitrary. By the splitting principle (cf.\ section (2.5) in
\cite{KoARR}), there is a morphism $p: Y' \ra Y$ such that 
$p^*: K_0(Y) \ra K_0(Y')$ and $p^*: K_0(\gS_l,Y) \ra K_0(\gS_l,Y')$
are finite free ring extensions and such that $p^*([\cF])$ 
decomposes into a sum of classes of invertible $\cO_{Y'}$-modules.
Since $\lambda_{-1}$ and $\tau^l$ commute with $p^*$ and since they
are multiplicative (see Proposition 1.7(a)), Corollary 1.11 now follows
from the case $\rank(\cF)=1$ proved above.

{\bf Remark 1.12}. Let $\gG$ be a group. We suppose that $\gG$ acts
on $X$ by $Y$-automorphisms. For any locally free $\cO_X$-module 
$\cF$ and for any $l\ge 0$, the $l$-th external tensor power 
$\cF^{\boxtimes l}$ then carries a natural action of the wreath
product $\gS_l\langle\gG\rangle := \gG^l \semidirect \gS_l$. Hence, as above,
one can construct operations 
\[\tau^l: K_q(\gG,X) \ra K_q(\gS_l\langle\gG\rangle,X^l), \quad 
q\ge 0, \quad l\ge 0.\]
One easily verifies that all statements proved so far hold in this
more general situation, too. Now let $\gG$ be the symmetric group
$\gS_k$ for some $k \ge 1$ and let $f: X^k \ra Y$ be the projection 
associated with a morphism $X\ra Y$ between noetherian schemes. Then,
induction with respect to the canonical embedding $\gS_l\langle
\gS_k\rangle \hookrightarrow \gS_{kl}$ yields homomorphisms
\[\Ind_{\gS_l\langle\gS_k\rangle}^{\gS_{kl}} : 
K_q(\gS_l\langle\gS_k\rangle,X^{kl}) \ra K_q(\gS_{kl}, X^{kl}),
\quad q\ge 0, \quad l\ge 0.\]
One easily verifies that the homomorphism $\prod_{l\ge 0} 
\Ind_{\gS_l\langle\gS_k\rangle}^{\gS_{kl}}$ commutes with cross 
products. Finally, the composition
\[K_q(\gS_k, X^k) \,\, \stackrel{\tau^l}{\longrightarrow}\,\,
K_q(\gS_l\langle\gS_k \rangle, (X^k)^l) \,\, 
\stackrel{\Ind_{\gS_l\langle \gS_k\rangle}^{\gS_{kl}}}{\longrightarrow}
\,\, K_q(\gS_{kl},X^{kl})\]
corresponds to the power operation mentioned in Remark 1.1.

Let $f= \id$ again. We now establish a fundamental relation between
the Adams operation $\psi^l$ on $K_q(Y)$ and the composition
\[K_q(Y) \,\, \stackrel{\tau^l}{\longrightarrow} \,\,
K_q(\gS_l, Y) \,\, \stackrel{\Res_{C_l}^{\gS_l}}{\longrightarrow}
\,\, K_q(C_l,Y)\]
which will be denoted by $\tau^l$ again. One should be able to prove
this relation also for the higher $K$-groups without the additional
assumptions made in the following proposition.

{\bf Proposition 1.13}. Let $l$ be a prime. Then the following diagram
commutes:
\[\begin{array}{ccc}
K_0(Y) & \stackrel{\tau^l}{\longrightarrow} & K_0(C_l, Y)\\
\\
{\ss \psi^l} \downarrow \phantom{\ss \psi^l} && 
\phantom{\ss {\rm can}} \downarrow {\ss {\rm can}}\\
\\
K_0(Y) & \stackrel{{\rm can}}{\longrightarrow} & 
K_0(C_l,Y)/[\cO_Y[I_l]] \cdot K_0(C_l,Y).
\end{array}\]
The corresponding diagrams of higher $K$-groups commute if $l$ is 
invertible on $Y$ and if $Y$ is affine or a smooth quasi-projective
scheme over an affine regular base $S$.

{\bf Proof}. First, we prove the assertion for Grothendieck groups.
For any $1\le i \le l-1$, the action of $C_l$ on $\gS_l/(\gS_i \times
\gS_{l-i})$ has no fixed points; thus, there is a $C_l$-stable 
system of representatives in $\gS_l$ for $\gS_l/(\gS_i \times
\gS_{l-i})$. This implies that, for any 
$(\gS_i \times \gS_{l-i})$-module $\cF$, the element 
$\left[\Ind_{\gS_i \times \gS_{l-i}}^{\gS_l} (\cF) \right]$ is
contained in the ideal of $K_0(C_l,Y)$ 
generated by $\left[\Ind_{\gS_i\times\gS_{l-i}}
^{\gS_l}(\cO_Y)\right]$ and then in the ideal generated by
$[\cO_Y[I_l]]$. It now follows from Corollaries 1.3 and 1.10 that
the composition
\[K_0(Y) \,\, \stackrel{\tau^l}{\longrightarrow} \,\,
K_0(C_l,Y) \,\, \stackrel{{\rm can}}{\longrightarrow} \,\,
K_0(C_l,Y)/[\cO_Y[I_l]]\cdot K_0(C_l,Y)\]
is a homomorphism. Since the same holds for the composition
\[K_0(Y) \,\, \stackrel{\psi^l}{\longrightarrow}\,\,
K_0(Y) \,\, \stackrel{{\rm can}}{\longrightarrow} \,\,
K_0(C_l,Y)/[\cO_Y[I_l]]\cdot K_0(C_l,Y),\]
it suffices to show the equality 
\[\tau^l([\cF]) = \psi^l([\cF]) \quad {\rm in} \quad
K_0(C_l,Y)/[\cO_Y[I_l]]\cdot K_0(C_l,Y)\] 
for any locally free 
$\cO_Y$-module $\cF$. We choose a morphism $p:Y' \ra Y$ for $\cF$
as in the proof of Proposition 1.11. Then, the homomorphism
\[p^*: K_0(C_l,Y)/[\cO_Y[I_l]]\cdot K_0(C_l,Y) \ra 
K_0(C_l, Y')/[\cO_{Y'}[I_l]]\cdot K_0(C_l,Y')\] is still
injective by (4.C)(ii) on p.\ 28 in \cite{Mat}. 
As in the proof of Corollary 1.11, it thus suffices 
to show the above equality for invertible $\cO_Y$-modules. In this
case, both sides are equal to $[\cF^{\otimes l}]$.\\
We now prove the assertion for higher $K$-groups. For this, let 
$\varepsilon$ denote the idempotent element $\varepsilon:= 
\frac{1}{l} \sum_{i=0}^l [c^i] \in \cO_Y[C_l]$, $G$ the group of 
automorphisms of $C_l$ and $C_l \semidirect G$ the semi-direct
product of $C_l$ with $G$ (defined by $(c^i,\gs) \cdot (c^j,\tau)
:= (c^{i+\gs j}, \gs \tau)$ for $i,j \in \ZZ/l\ZZ$ and 
$\gs, \tau \in (\ZZ/l\ZZ)^\times \cong G$). For any 
$(C_l \semidirect G)$-module $\cF$ on $Y$ with $\varepsilon \cF = 0$,
the module $\cF^G$ of $G$-fixed elements will be considered as 
a $C_l$-module with trivial $C_l$-action. Then, obviously, the map
\[\alpha: \cF \ra \cF^G \otimes \cO_Y[C_l], \quad 
f \mapsto \sum_{i=0}^{l-1}\left(\sum_{\gs \in G} 
\gs(c^{-i}(f))\right) \otimes [c^i],\]
is a homomorphism of $C_l$-modules on $Y$ and the sequence
\[0 \ra \cF \,\, \stackrel{\alpha}{\ra}\,\, 
\cF^G \otimes \cO_Y[C_l] \,\, \stackrel{\gS}{\ra} \,\, 
\cF^G \ra 0 \]
is a split short sequence of $C_l$-modules. We now consider the
functors
\[\begin{array}{cccccc}
& \cP(C_l \semidirect G, Y) & \ra & \cP(Y) & 
\stackrel{{\rm can}}{\ra} & \cP(C_l, Y) \\
F_0: & \cF & \mapsto & \varepsilon \cF \\
F_1: & \cF & \mapsto & ((1-\varepsilon)\cF)^G.
\end{array}\]
Then, the direct sum $F_0 \oplus F_1 \otimes \cH_{Y,l}$ is isomorphic
to the forgetful functor 
\[\cP(C_l \semidirect G, Y) \ra \cP(C_l,Y).\]
Hence, the functors $F_0$ and $F_1$ are exact and induce homomorphisms
\[F_0 \quad {\rm and} \quad F_1: K_q(C_l\semidirect G, Y) \doublearrow
K_q(Y)\]
with $F_0(y) + [\cH_{Y,l}] \cdot F_1(y) = y$ in $K_q(C_l,Y)$ for all
$y \in K_q(C_l\semidirect G,Y)$. In particular, the diagram 
\[\begin{array}{ccc}
K_q(C_l\semidirect G, Y) & 
\stackrel{{\rm Res}^{C_l \rtimes G}_{C_l}}{\longrightarrow} &
K_q(C_l, Y) \\
\\
{\ss F_0 - F_1} \downarrow \phantom{\ss F_0 - F_1} && 
\phantom{\ss {\rm can}} \downarrow {\ss {\rm can}}\\
\\
K_q(Y) & \stackrel{{\rm can}}{\longrightarrow} & 
K_q(C_l, Y)/[\cO_Y[I_l]] \cdot K_q(C_l,Y)
\end{array}\]
commutes.
We now consider $C_l \semidirect G$ as a subgroup of $\gS_l$ in the
obvious way.
Then it is easy to verify (see also Remark 1.6) that the compositions
\[K_q(Y) \,\, \stackrel{\tau^l}{\longrightarrow} \,\,
K_q(\gS_l, Y) \,\, \stackrel{\Res^{\gS_l}_{C_l\semidirect G}}
{\longrightarrow} \,\, K_q(C_l\semidirect G, Y) \,\,
\stackrel{F_0, F_1}{\longrightarrow} \,\, K_q(Y)\]
are equal to the cyclic power operations $[0]_l$ and $[1]_l$
constructed in section 3 of \cite{KoAdPro}. By Corollary 2.2(c) on
p.\ 142 in \cite{KoAdHi} (applied to the case $\gG =\{1\}$), the 
difference between these compositions equals $\psi^l$, if $Y$ is affine.
By Jouanolou's trick (e.g., see Satz (4.4) on p.\ 211 in 
\cite{KoARR}) and Quillen's resolution theorem (in the form of Satz (2.1)
on p.\ 195 in \cite{KoARR}), the same holds if $Y$ is a 
smooth quasi-projective
scheme over an affine regular base $S$. This proves Proposition 1.13
for higher $K$-groups.

Now let $k$ be an algebraically closed field of characteristic $0$ and
$l$ be a prime. In the subsequent proposition, we will explicitly 
describe the $l$-th tensor power operation
\[\tau^l: K_1(k) \,\, \stackrel{\tau^l}{\longrightarrow}\,\,
K_1(\gS_l, k) \,\, \stackrel{\Res^{\gS_l}_{C_l}}{\longrightarrow}
\,\, K_1(C_l, k).\]
For this, we fix a non-trivial character $\chi: C_l \ra k^\times$. 
For any $i=0, \ldots, l-1$, let $V_i$ denote the one-dimensional
representation of $C_l$ given by the character $\chi^i$. Then it 
is well-known that the homomorphism
\[\begin{array}{ccc}
(k^\times)^{\ZZ/l\ZZ} & \ra & K_1(k[C_l]) = K_1(C_l, k)\\
(\beta_0, \ldots, \beta_{l-1}) & \mapsto & 
(V_0, \beta_0) \oplus \ldots \oplus (V_{l-1}, \beta_{l-1})
\end{array}\]
is bijective. Here, as usual, the pairs $(V,\beta)$, $V$ a finite
dimensional $k[C_l]$-module and $\beta \in \Aut_{k[C_l]}(V)$,
are considered as generators of the group $K_1(k[C_l])$. 

{\bf Proposition 1.14}. Let $l$ be a prime. Then the following diagram
commutes:
\[\begin{array}{ccc}
K_1(k) & \stackrel{\tau^l}{\ra} & K_1(C_l, k)\\
\\
\| && \| \\
\\
k^\times & \ra & (k^\times)^{\ZZ/l\ZZ}\\
\beta & \mapsto & (\beta^{l-1}, \beta^{-1}, \ldots, \beta^{-1}).
\end{array}\]

{\bf Proof}. Let $K_0(\ZZ, k[C_l])$ denote the Grothendieck group of 
all pairs $(V, \beta)$ as above. An obvious generalization of the
construction of this section (see also Remark 1.12) yields a tensor
power operation
\[\tau^l: K_0(\ZZ, k) \ra K_0(\ZZ, k[C_l])\]
such that, for any vector space $V$ over $k$ and any $\beta \in
\Aut_k(V)$, the class $[(V,\beta)]$ is mapped to the class 
$[(V^{\otimes l}, \beta^{\otimes l})]$ where $V^{\otimes l}$ is
considered as a $C_l$-module as usual. By restricting, we obtain
a tensor power operation
\[\tau^l: \tilde{K}_0(\ZZ,k) \ra \tilde{K}_0(\ZZ, k[C_l])\]
between the reduced Grothendieck groups $\tilde{K}_0(\ZZ,k) :=
\ker(K_0(\ZZ,k) \,\, \stackrel{{\rm can}}{\ra} \,\, K_0(k))$
and $\tilde{K}_0(\ZZ, k[C_l]) := \ker(K_0(\ZZ,k[C_l]) \,\,
\stackrel{{\rm can}}{\ra} \,\, K_0(k[C_l]))$. As in Theorem 3.3 on
p.\ 145 in \cite{KoAdHi} one shows that the following diagram
commutes:
\[\begin{array}{ccc}
\tilde{K}_0(\ZZ,k) & \stackrel{{\rm can}}{\longrightarrow} &
K_1(k)\\
\\
{\ss \tau^l} \downarrow \phantom{\ss \tau^l} && 
\phantom{\ss \tau^l} \downarrow {\ss \tau^l} \\
\\
\tilde{K}_0(\ZZ,k[C_l]) & \stackrel{{\rm can}}{\longrightarrow} &
K_1(k[C_l]).
\end{array}\]
By Corollary 1.10, we have the following equality in 
$\tilde{K}_0(\ZZ,k[C_l])$ for all $\beta \in k^\times$:
\begin{eqnarray*}
\lefteqn{\tau^l\left([(k,\beta)]-[(k,1)]\right) = }\\
&& = \sum_{i=0}^l (-1)^{l-i}\left[\Res^{\gS_l}_{C_l} \left(
\Ind_{\gS_i \times \gS_{l-i}}^{\gS_l}(k\otimes k_\sgn, \beta^i)\right)
\right] \\
&& = [(V_0, \beta^l)] + \sum_{i=1}^{l-1} (-1)^{l-i} \frac{1}{l}
{l \choose i}[(k[C_l],\beta^i)] + (-1)^l[(V_{0,\sgn},1)].
\end{eqnarray*}
(The last equality follows from the fact that, for all $i=1, \ldots,
l-1$, there is a $C_l$-stable
system of representatives in $\gS_l$ for $\gS_l/(\gS_i \times
\gS_{l-i})$.) Because of
\[\sum_{i=1}^{l-1} (-1)^{l-i} \frac{1}{l}{l \choose i}i = 
\sum_{i=1}^{l-1} (-1)^{l-i} \frac{(l-1)!}{(i-1)!(l-i)!} =
\sum_{i=0}^{l-2}(-1)^{l-i-1} {l-1 \choose i} = -1\]
we have 
\[\tau^l((k,\beta)) = 
\tau^l((k,\beta)-(k,1)) = (V_0, \beta^l) - (k[C_l], \beta) \quad
{\rm in} \quad K_1(C_l, k).\]
This proves Proposition 1.14.

\bigskip

\section*{\S 2 A K\"unneth Formula}

Let $l \in \NN$, and let $Y$ be a noetherian scheme on which each coherent
module is a quotient of a locally free $\cO_Y$-module. Let $f:X \ra Y$
be a flat projective local complete intersection morphism (cf.\ p.\ 86
in \cite{FL}). Being the composition of the successive projections
$X^i \ra X^{i-1}$, $i=1, \ldots, l$, the morphism
\[f^l: X^l \ra Y\]
is then of the same kind (see Proposition 3.12 on p.\ 87 in
\cite{FL}). 
By section (2.6) in \cite{KoARR}, we have push-forward homomorphisms
\[f_*: K_q(X) \ra K_q(Y) \quad {\rm and} \quad 
f^l_*:  K_q(\gS_l, X^l) \ra K_q(\gS_l,Y)\]
(for all $q \ge 0$) whose definition will be recalled in the
subsequent proofs. The aim of this section is to prove the following
theorem:

{\bf Theorem 2.1} (K\"unneth formula). For any $q \ge 0$, the
following diagram commutes:
\[\begin{array}{ccc}
K_q(X) & \stackrel{\tau^l}{\longrightarrow} & K_q(\gS_l, X^l)\\
\\
{\ss f_*} \downarrow \phantom{\ss f_*} && 
\phantom{\ss f^l_*} \downarrow {\ss f^l_*}\\
\\
K_q(Y) & \stackrel{\tau^l}{\longrightarrow} & K_q(\gS_l, Y).
\end{array}\]

{\bf Proof}. There is a factorization $X \,\, 
\stackrel{i}{\hookrightarrow} \,\, \PP_Y(\cE) \,\,
\stackrel{p}{\ra} \,\, Y$ of $f$ into a regular closed immersion $i$
and the projection of the projective space bundle associated with
a locally free $\cO_Y$-module $\cE$. Since we have $f_* = p_* \circ
i_*$ and $f^l_* = p^l_* \circ i^l_*$, 
Theorem 2.1 immediately follows from the subsequent 
Propositions 2.2 and 2.3. As in the proof of Proposition 1.4, we 
will use the notations of \cite{KoSh} in their proof.

{\bf Proposition 2.2}. Let $i: X \hookrightarrow \tilde{X}$ be a 
regular closed immersion of flat quasi-projective schemes over $Y$.
Then the following diagram commutes for all $q\ge 0$:
\[\begin{array}{ccc}
K_q(X) & \stackrel{\tau^l}{\longrightarrow} & K_q(\gS_l, X^l)\\
\\
{\ss i_*} \downarrow \phantom{\ss i_*} && 
\phantom{\ss i^l_*} \downarrow {\ss i^l_*}\\
\\
K_q(\tilde{X}) & \stackrel{\tau^l}{\longrightarrow} &
K_q(\gS_l, \tilde{X}^l).
\end{array}\]

{\bf Proof}. Being the composition of the regular closed immersions
\[X^l=X\times \ldots \times X \times X \hookrightarrow 
X\times \ldots \times X \times \tilde{X}, \quad \ldots,\qquad \qquad \qquad\] 
\[\qquad \qquad \qquad 
X\times \tilde{X} \times \ldots \times \tilde{X} \hookrightarrow
\tilde{X}\times \tilde{X} \times \ldots \times \tilde{X} = 
\tilde{X}^l,\]
the morphism $i^l: X^l \ra \tilde{X}^l$ is a regular closed immersion,
too. Let $\cP_\infty(\tilde{X})$ (respectively, $\cP_{\infty, {\rm fl}}
(\tilde{X})$) denote the exact category of all coherent 
$\cO_{\tilde{X}}$-modules which possess a finite resolution by locally
free $\cO_{\tilde{X}}$-modules (respectively, 
which are in addition flat over $Y$). 
Similarly, let $\cP_\infty(\gS_l, \tilde{X}^l)$ and
$\cP_{\infty, {\rm fl}}(\gS_l, \tilde{X}^l))$ denote the exact
categories of all coherent $\gS_l$-modules on $\tilde{X}^l$ which,
after forgetting the $\gS_l$-structure, are contained in  
$\cP_\infty(\tilde{X}^l)$ and 
$\cP_{\infty, {\rm fl}}(\tilde{X}^l)$,
respectively. By Lemma (2.2) in \cite{KoARR} and Lemma (3.4)(b) in
\cite{KoGRR}, all modules in $\cP_\infty(\gS_l, \tilde{X}^l)$ then
possess a finite resolution by locally free $\gS_l$-modules  on 
$\tilde{X}^l$. The inclusions
\[\cP(\tilde{X}) \hookrightarrow \cP_{\infty, {\rm fl}}(\tilde{X})
\hookrightarrow \cP_\infty(\tilde{X}) \quad
{\rm and}
\quad \cP(\gS_l,\tilde{X}^l) \hookrightarrow 
\cP_{\infty, {\rm fl}}(\gS_l, \tilde{X}^l) \hookrightarrow
\cP_\infty(\gS_l,\tilde{X}^l)\]
induce homotopy equivalences between the geometric realizations 
of the corresponding ($l$-fold iterated) $G$-constructions by 
Quillen's resolution theorem (see Corollary 1 on
p.\ 109 in \cite{Q}). Furthermore, since $X$ is flat over $Y$, 
we have well-defined exact functors
\[i_*: \cP(X) \ra \cP_{\infty, {\rm fl}}(\tilde{X}) \quad {\rm and}
\quad i^l_*: \cP(\gS_l,X^l) \ra \cP_{\infty, {\rm fl}}(\gS_l, 
\tilde{X}^l)\]
(see section (2.6) in \cite{KoARR}). By definition, the maps
\[i_*: K_q(X) \ra K_q(\tilde{X}) \quad {\rm and} \quad
i^l_*: K_q(\gS_l,{X}^l) \ra K_q(\gS_l,\tilde{X}^l)\]
are induced by the compositions
\[|G\cP(X)| \,\, \stackrel{i_*}{\ra} \,\,
|G\cP_\infty(\tilde{X})| \simeq |G\cP(\tilde{X})| \]
and
\[|G^l\cP(\gS_l,X^l) \,\, 
\stackrel{i^l_*}{\ra} \,\, 
|G^l\cP_\infty(\gS_l,\tilde{X}^l)| \simeq 
|G^l\cP(\gS_l, \tilde{X}^l)|.\]
For any $i,j \in \NN$ and any 
$\cE_1 \in \cP_{\infty, {\rm fl}}(\tilde{X}^i)$, 
$\cE_2 \in \cP_{\infty, {\rm fl}}(\tilde{X}^j)$, the external
tensor product $\cE_1 \boxtimes \cE_2$ is obviously flat over $Y$,
again. Furthermore, we have: If $\cF^\cdot_i \ra \cE_i$, $i=1,2$, are
finite locally free resolutions, then the total complex associated
with the double complex $\cF^\cdot_1 \boxtimes \cF^\cdot_2$ is a finite 
locally free resolution of $\cE_1 \boxtimes \cE_2$ since $\cE_1$ and
$\cE_2$ are flat over $Y$ and since $\tilde{X}$ is flat over $Y$.
Finally, the functor
\[\cP_{\infty,{\rm fl}}(\tilde{X}^i) \times
\cP_{\infty, {\rm fl}}(\tilde{X}^j) \ra 
\cP_{\infty, {\rm fl}}(\tilde{X}^{i+j}), \quad
(\cE_1, \cE_2) \mapsto \cE_1 \boxtimes \cE_2,\]
is obviously bi-exact. This implies that we have functors 
\[\times: \cP_{\infty, {\rm fl}}(\gS_i, \tilde{X}^i) \times
\cP_{\infty, {\rm fl}}(\gS_j, \tilde{X}^j) \ra 
\cP_{\infty, {\rm fl}}(\gS_{i+j}, \tilde{X}^{i+j}) \]
and
\[F_k(\cP_{\infty, {\rm fl}}(\tilde{X})) \ra 
\cP_{\infty, {\rm fl}}(\gS_k, \tilde{X}^k)\]
as in section 1 which satisfy axioms (E1) through (E5) in
section 7 of \cite{GrEx}. Now, Grayson's construction in section 7
of \cite{GrEx} yields a simplicial map $\tau^l: 
G\cP_{\infty,{\rm fl}}(\tilde{X}) \ra \cH^l\cP_{\infty,{\rm fl}}
(\gS_l,\tilde{X}^l)$ and then an $l$-th external tensor power
operation
\[|G\cP_{\infty, {\rm fl}}(\tilde{X})| \,\,
\stackrel{|\tau^l|}{\ra} \,\,
|\cH^l\cP_{\infty,{\rm fl}}(\gS_l, \tilde{X}^l)| \,\,
\stackrel{\Xi^l}{\ra}\,\,
|G^l\cP_{\infty,{\rm fl}}(\gS_l, \tilde{X}^l)|\]
for $\cP_{\infty,{\rm fl}}(\tilde{X})$ which is compatible with the
corresponding operation for $\cP(\tilde{X})$ defined in section 1.
Finally, for any $\cE_1, \ldots, \cE_l \in \cP(X)$, the 
canonical homomorphism
\[i_*(\cE_1) \boxtimes \ldots \boxtimes i_*(\cE_l) \ra 
i^l_*(\cE_1 \boxtimes \ldots \boxtimes \cE_l)\]
is bijective as one can easily check stalk-wise. Altogether, we have
the following (up to homotopy) commutative diagram of continuous maps
in which the left and right lower arrow are homotopy equivalences:
\[\begin{array}{ccccc}
|G\cP(X)| & \stackrel{\tau^l}{\ra} &
|\cH^l\cP(\gS_l,X^l)| & \stackrel{\Xi^l}{\ra} &
|G^l\cP(\gS_l, X^l)|\\
\\
{\ss i_*} \downarrow \phantom{\ss i_*} && 
{\ss i^l_*} \downarrow \phantom{\ss i^l_*}&&
\phantom{\ss i^l_*}\downarrow {\ss i^l_*} \\
\\
|G^l\cP_{\infty,{\rm fl}}(\tilde{X})| &
\stackrel{\tau^l}{\ra} & 
|\cH^l\cP_{\infty, {\rm fl}}(\gS_l, \tilde{X}^l| & 
\stackrel{\Xi^l}{\ra}&
|G^l\cP_{\infty,{\rm fl}}(\gS_l,\tilde{X}^l)|\\
\\
\uparrow && \uparrow && \uparrow\\
\\
|G\cP(\tilde{X})| & \stackrel{\tau^l}{\ra} &
|\cH^l\cP(\gS_l, \tilde{X}^l)| & \stackrel{\Xi^l}{\ra}&
|G^l\cP(\gS_l, \tilde{X}^l)|.
\end{array}\]
This proves Proposition 2.2.

{\bf Proposition 2.3}. Let $p: \PP:= \PP_Y(\cE) \ra Y$ be the
projective space bundle associated with a locally free $\cO_Y$-module
$\cE$. Then the following diagram commutes for all $q\ge 0$:
\[\begin{array}{ccc}
K_q(\PP)&\stackrel{\tau^l}{\longrightarrow} & K_q(\gS_l,\PP^l)\\
\\
{\ss p_*} \downarrow \phantom{\ss p_*} && 
\phantom{\ss p^l_*} \downarrow {\ss p^l_*}\\
\\
K_q(Y) & \stackrel{\tau^l}{\longrightarrow}& K_q(\gS_l,Y).
\end{array}\]

{\bf Proof}. Let $\cR(\PP)$ denote the exact category consisting of
all locally free $\cO_\PP$-modules $\cF$ with $R^jp_*\cF(i)=0$ for 
all $j\ge 1$ and $i \ge 0$. By Lemma 1.13 on p.\ 141 in \cite{Q},
we have a well-defined exact functor $p_*: \cR(\PP) \ra \cP(Y)$.
The inclusion $\cR(\PP) \hookrightarrow \cP(\PP)$ induces a homotopy
equivalence $|G\cR(\PP)| \ra |G\cP(\PP)|$ by Lemma 2.2 on p.\ 142
in \cite{Q}. By definition, the map $p_*: K_q(\PP) \ra K_q(Y)$ is
induced by the composition $|G\cP(\PP)| \simeq |G\cR(\PP)| \,\,
\stackrel{p_*}{\ra}\,\, |G\cP(Y)|$. \\
Now, let $\cR(\gS_l, \PP^l)$ denote the exact category consisting of
all locally free $\gS_l$-modules $\cF$ on $\PP^l$ for which 
$R^jp^l_*\cF(i)$ vanishes for all $j \ge 1$ and $i \ge 0$ and for 
which $p^l_*\cF(i)$ is a locally free $\cO_Y$-module for all $i\ge 0$
(for this, we put $\cF(i) := \cF \otimes \cO(i)^{\boxtimes l}$). In 
particular, the functor $p^l_*: \cR(\gS_l, \PP^l) \ra \cP(\gS_l, Y)$
is well-defined and exact. For any $\cF_1, \ldots, \cF_l \in
\cR(\PP)$, the higher direct image $R^jp^l_*(\cF_1 \boxtimes
\ldots \boxtimes \cF_l)(i)$ vanishes for all $j\ge 1$ and $i\ge 0$,
and the module $p^l_*(\cF_1 \boxtimes \ldots \boxtimes \cF_l)(i)
\cong p_*\cF_1(i) \otimes \ldots \otimes p_*\cF_l(i)$ is locally
free for all $i\ge 0$ (by Proposition 9.3 on p.\ 255 in \cite{H}).
Similarly to the proof of Proposition 2.2, we can therefore define
a simplicial map $\tau^l: G\cR(\PP) \ra \cH^l\cR(\gS_l,\PP^l)$ and
an $l$-th external tensor power operation
\[|G\cR(\PP)| \,\, \stackrel{|\tau^l|}{\ra} \,\,
|\cH^l\cR(\gS_l, \PP^l)| \,\,
\stackrel{\Xi^l}{\ra}\,\, |G^l\cR(\gS_l,\PP^l)|\]
such that the following diagram commutes (up to homotopy):
\[\begin{array}{ccccc}
|G\cP(\PP)| & \stackrel{|\tau^l|}{\ra} &
|\cH^l\cP(\gS_l, \PP^l)| & \stackrel{\Xi^l}{\ra}&
|G^l\cP(\gS_l,\PP^l)|\\
\\
\uparrow && \uparrow && \uparrow\\
\\
|G\cR(\PP)| & \stackrel{|\tau^l|}{\ra}&
|\cH^l\cR(\gS_l,\PP^l)| & \stackrel{\Xi^l}{\ra}&
|G^l\cR(\gS_l,\PP^l)|\\
\\
{\ss p_*} \downarrow \phantom{\ss p_*} && 
{\ss p^l_*} \downarrow \phantom{\ss p^l_*} &&
\phantom{\ss p^l_*} \downarrow {\ss p^l_*}\\
\\
|G\cP(Y)| & \stackrel{|\tau^l|}{\ra} & 
|\cH^l\cP(\gS_l,Y)| & \stackrel{\Xi^l}{\ra}& 
|G^l\cP(\gS_l,Y)|.
\end{array}\]

Now, let $r:\RR := \PP_Y(\cE^{\otimes l}) \ra Y$ denote the projective
space bundle associated with $\cE^{\otimes l}$; it will be considered 
as a $\gS_l$-scheme over $Y$ by virtue of the usual action of $\gS_l$
on $\cE^{\otimes l}$. The canonical $\gS_l$-epimorphism 
$(p^l)^*(\cE^{\otimes l}) \rightepi \cO(1)^{\boxtimes l}$ defines a regular
closed $\gS_l$-immersion $i: \PP^l \ra \RR$ (see section (1.3) in
\cite{KoARR} and Proposition 3.11 on p.\ 85 in \cite{FL}). Let 
$\cR\cP_{\infty, {\rm fl}}(\gS_l,\RR)$ denote the exact category
consisting of all coherent $\gS_l$-modules $\cF$ on $\RR$ which
possess a finite resolution by locally free $\cO_\RR$-modules, which
are flat over $Y$, and for which $R^jr_*\cF(i)$ vanishes for all 
$j\ge 1$ and $i \ge 0$. As in the proof of Proposition 2.2, we then
have a well-defined exact functor $i_*: \cR(\gS_l,\PP^l) \ra 
\cR\cP_{\infty, {\rm fl}}(\gS_l, \RR)$. We now prove that the 
inclusion $\cR\cP_{\infty, {\rm fl}}(\gS_l, \RR) \hookrightarrow
\cP_{\infty, {\rm fl}}(\gS_l,\RR)$ induces a homotopy equivalence
\[|G^l\cR\cP_{\infty, {\rm fl}}(\gS_l, \RR)| \ra 
|G^l\cP_{\infty, {\rm fl}}(\gS_l, \RR)|\]
(here, $\cP_{\infty, {\rm fl}}(\gS_l,\RR)$ is defined as in the proof
of Proposition 2.2). For this, let $\cR_n$ denote the exact category
consisting of all $\cF \in \cP_{\infty, {\rm fl}}(\gS_l, \RR)$ with
$\cF(n) \in \cR\cP_{\infty, {\rm fl}}(\gS_l, \RR)$ (for any $n \ge 0$). 
Thus we have $\cR_0 = \cR\cP_{\infty, {\rm fl}}(\gS_l,\RR)$ and
$\cR_0 \subseteq \cR_1 \subseteq \cR_2 \subseteq \ldots$. By Lemma
1.12 on p.\ 141 in \cite{Q}, we have $\cup_{n \ge 0} \cR_n =
\cP_{\infty, {\rm fl}}(\gS_l, \RR)$. (Note that
Lemma 1.12 in \cite{Q} holds even for coherent $Y$-flat modules 
$\cF$. Quillen's proof for locally free $\cF$ need not to be 
changed.) By property (9) on p.\ 104 in \cite{Q},
it therefore suffices to show that the inclusion $\cR_n
\hookrightarrow \cR_{n+1}$ induces a homotopy equivalence 
$|G^l\cR_n| \ra |G^l\cR_{n+1}|$ for all $n\ge 0$. The exact
functors
\[u_p: \cR_{n+1} \ra \cR_n, \quad 
\cF \mapsto r^* \gL^p((\cE^{\otimes l})^\vee) \otimes \cF(p), 
\quad p\ge 1,\]
induce homomorphisms $u_p: K_q(\cR_{n+1}) \ra K_q(\cR_n)$, $p\ge 1$,
(for all $q \ge 0$). The functorial exact sequence
\[0 \ra \cF \ra r^*(\cE^{\otimes l})^\vee \otimes \cF(1) \ra \ldots
\ra \cF(\rank(\cE^{\otimes l})) \ra 0 \]
(see sequence (2.2) on p.\ 107 in \cite{FL}) 
together with Corollary 3 of Theorem 2 on
p.\ 107 in \cite{Q} imply  that $\sum_{p\ge 1} (-1)^{p-1}u_p$
is inverse to the canonical map $K_q(\cR_n) \ra 
K_q(\cR_{n+1})$. Hence, the map $|G^l\cR_n| \ra |G^l\cR_{n+1}|$
is a homotopy equivalence. It finally follows from Lemma 1.13 on
p.\ 141 in \cite{Q} (which again holds more generally for
coherent $Y$-flat modules) that we have well-defined exact functors
$r_*: \cR\cP_{\infty,{\rm fl}}(\gS_l,\RR) \ra \cP(\gS_l,Y)$ and
$r_*: \cR(\gS_l, \RR) \ra \cP(\gS_l,Y)$. \\
Altogether, we now have the following (up to natural equivalence)
commutative diagram of exact categories and exact functors in which
all inclusions in the upper right corner become homotopy equivalences
after passing to the geometric realization of the ($l$-fold iterated)
$G$-construction:
\[\begin{array}{ccccc}
\cP(\gS_l, \PP^l) & \stackrel{i_*}{\ra} &
\cP_\infty(\gS_l,\RR)\\
\\
&& \uparrow & \nwarrow\\
\\
\uparrow && \cP_{\infty, {\rm fl}}(\gS_l, \RR) & 
\leftarrow & \cP(\gS_l, \RR)\\
\\
&& \uparrow & & \uparrow\\
\\
\cR(\gS_l,\PP^l) &\stackrel{i_*}{\ra} & \cR\cP_{\infty, {\rm fl}}
(\gS_l, \RR) & \leftarrow & \cR(\gS_l, \RR)\\
\\
& {\ss p^l_*} \searrow \phantom{\ss p^l_*} & 
{\ss r_*} \downarrow \phantom{\ss r_*} &
\phantom{\ss r_*} \swarrow {\ss r_*}\\
\\
&& \cP(\gS_l,Y).
\end{array}\]
This diagram together with the diagram above prove Proposition 2.3.

{\bf Remark 2.4}. The following generalizations of Theorem 2.1 can
be proved with the same methods:\\
(a) Let $q\ge 0$. Let $K_q^{Y-{\rm fl}}(X)$ denote the $q$-th
$K$-group associated with the exact category consisting of all
coherent $Y$-flat $\cO_X$-modules. Then one can naturally define 
a push-forward homomorphism $f_*: K_q^{Y-{\rm fl}}(X) \ra K_q(Y)$ and
an $l$-th external tensor power operation $\tau^l: 
K_q^{Y-{\rm fl}}(X) \ra K_q^{Y-{\rm fl}}(\gS_l, X^l)$ such that the
following diagram commutes:
\[\begin{array}{ccc}
K_q^{Y-{\rm fl}}(X) & \stackrel{\tau^l}{\longrightarrow} &
K_q^{Y-{\rm fl}}(\gS_l, X^l) \\
\\
{\ss f_*} \downarrow \phantom{\ss f_*} && 
\phantom{\ss f^l_*} \downarrow {\ss f^l_*}\\
\\
K_q(Y) & \stackrel{\tau^l}{\longrightarrow} & K_q(\gS_l,Y).
\end{array}\]
(b) Let $Y$ be a noetherian scheme and $X_1, X_2$ flat noetherian
schemes over $Y$ on which each coherent module is a quotient of a
locally free module. Furthermore, let $f:X_2 \ra X_1$ be 
a projective local complete intersection morphism. Then the following 
diagram commutes for all $q \ge 0$:
\[\begin{array}{ccc}
K_q(X_2) & \stackrel{\tau^l}{\longrightarrow} & K_q(\gS_l,X_2^l)\\
\\
{\ss f_*} \downarrow \phantom{\ss f_*} && 
\phantom{\ss f^l_*} \downarrow {\ss f^l_*}\\
\\
K_q(X_1) & \stackrel{\tau^l}{\longrightarrow} & K_q(\gS_l, X_1^l).
\end{array}\]
(Here, the fibred products $X_1^l$ and $X_2^l$ are formed over $Y$.)

\bigskip

\section*{\S 3 On a Certain Multiplicative Class in $K_0(\gS_l, X)$}

Let $l\in \NN$, and let $f: X \ra Y$ be a smooth quasi-projective morphism
between noetherian schemes. Let $\gO_{X/Y}$ denote the locally free
sheaf of relative differentials, $\gD: X \ra X^l$ the diagonal, and
$I$ the ideal of $K_0(\gS_l, Y)$ generated by the elements 
$[\gL^i(\cO_Y[I_l])]$, $i=1, \ldots, l-1$. The aim of this section is
to prove the following theorem:

{\bf Theorem 3.1}. The element $\lambda_{-1}([\gO_{X/Y} \otimes
\HXl]) \in K_0(\gS_l, X)$ is invertible in 
\[K_0(\gS_l,X)[l^{-1}]/IK_0(\gS_l,X)[l^{-1}]\] 
and in particular in
$K_0(C_l,X)[l^{-1}]/IK_0(C_l,X)[l^{-1}]$, 
and we have:
\[\gD_*\left(\lambda_{-1}([\gO_{X/Y}\otimes \HXl])^{-1}\right)
= 1 \quad {\rm in}
\quad K_0(C_l, X^l)[l^{-1}]/IK_0(C_l,X^l)[l^{-1}].\]

In subsections 1 through 5,  we prove some auxiliary results. The 
proof of Theorem 3.1 will be contained in subsection 6.

\subsection*{1. The Bott element $\gt^l(\cF)$ and the class
$\lambda_{-1}([\cF\otimes \HXl])$}

Let $l\in \NN$ and $X$ a $\gS_l$-scheme. For any locally free 
$\gS_l$-module $\cF$ on $X$, let $\gt^l(\cF) \in K_0(\gS_l, X)$
denote the $l$-th Bott element associated with $\cF$ (e.g., 
see section 4 in \cite{KoGRR}). It is characterized by the following 
properties: If $\rank(\cF) = 1$, then 
\[\gt^l(\cF) = [\cO_X] + [\cF] + \ldots + [\cF^{\otimes (l-1)}].\] 
If $\cF'$ is another
locally free $\gS_l$-module on $X$, then $\gt^l(\cF \oplus \cF')
= \gt^l(\cF) \cdot \gt^l(\cF')$. Finally, we have $\gt^l(p^*\cF) =
p^*\gt^l(\cF)$ in $K_0(\gS_l, X')$ for any morphism $p:X' \ra X$
between $\gS_l$-schemes. 

If $l$ is prime and $C_l \subseteq \gS_l$ acts on $X$ 
trivially, we have in $K_0(C_l,X)/([\cO_X[C_l]])$:
\begin{eqnarray*}
\lefteqn{\gt^l(\cF) \cdot \lambda_{-1}([\cF]) = 
\psi^l(\lambda_{-1}([\cF]))}\\
&&= \tau^l(\lambda_{-1}([\cF]))\\
&&= \lambda_{-1}([\cF\otimes \HXl]) \cdot \lambda_{-1}([\cF]).
\end{eqnarray*}
Here, we have used the splitting principle (see section (2.5) in
\cite{KoARR}), Proposition 1.13 and Corollary 1.11, successively.
The following proposition streng\-thens this computation:

{\bf Proposition 3.2}. For any locally free $\gS_l$-module $\cF$
on $X$, we have:
\[\gt^l(\cF) = \lambda_{-1}([\cF \otimes \HXl]) \quad {\rm in}
\quad K_0(\gS_l,X)/\left([\gL^i(\cO_X[I_l])], i=1, \ldots, l\right).\]

{\bf Proof}. For any $i=0, \ldots l-1$, we have
\begin{eqnarray*}
\lefteqn{[\gL^i(\HXl)]= \lambda^i([\cO_X[I_l]]-1)}\\
&&= \sum_{j=0}^i (-1)^j\left[\gL^{i-j}(\cO_X[I_l])\right]\\
&& =(-1)^i \quad {\rm in} \quad K_0(\gS_l,X)/\left([\gL^i(\cO_X[I_l])], i=1,
\ldots, l-1\right)
\end{eqnarray*}
since $\lambda^j(-1)=(-1)^j$ for $j\ge 0$.
Hence, we have for all invertible $\gS_l$-modules $\cF$ on $X$:
\begin{eqnarray*}
\lefteqn{\gt^l(\cF) = \sum_{i=0}^{l-1} [\cF^i]}\\
&&= \sum_{i=0}^{l-1} (-1)^i [\cF^i \otimes \gL^i(\HXl)]\\
&&= \sum_{i=0}^{l-1} (-1)^i [\gL^i(\cF\otimes \HXl)]\\
&&= \lambda_{-1}([\cF \otimes \HXl]) \quad {\rm in} \quad
K_0(\gS_l,X)/\left([\gL^i(\cO_X[I_l])], i=1, \ldots, l-1\right).
\end{eqnarray*}
As in Proposition 1.13, Proposition 3.2 now follows from this using
the splitting principle. 

\subsection*{2. On the cohomology of Koszul complexes}

Let $G$ be a group, and let $i: X \hookrightarrow \tilde{X}$ be a regular
closed $G$-immersion of noetherian $G$-schemes. Let $\cI := 
\ker(\cO_{\tilde{X}} \ra i_*(\cO_X))$ denote the corresponding
$G$-stable ideal in $\cO_{\tilde{X}}$. Furthermore, let $\cF$ be a 
locally free $G$-module on $X$ of rank $n$ and $\varepsilon: \cF \rightepi
\cI$ an epimorphism of $G$-modules. For any $i=0, \ldots, n$, let
$\cH^i(\gL^\cdot(\cF), d_\cdot)$ denote the $i$-th homology module
of the Koszul complex
\[0 \ra \gL^n(\cF) \,\, \stackrel{d_n}{\ra} \,\,
\ldots \,\, \stackrel{d_3}{\ra} \,\, \gL^2(\cF) \,\,
\stackrel{d_2}{\ra} \,\, \cF \,\, \stackrel{d_1=\varepsilon}{\ra}
\,\, \cO_{\tilde{X}} \ra  0\]
associated with the homomorphism $\varepsilon: \cF \,\,\rightepi\,\, \cI 
\subseteq \cO_{\tilde{X}}$. Furthermore, let the locally free
$G$-module $\cE$ on $X$ be defined by the short exact sequence
\[0 \ra \cE \ra i^*(\cF) 
\stackrel{i^*(\varepsilon)}{\ra} i^*(\cI) = \cI/\cI^2 \ra 0.\]

{\bf Lemma 3.3}. For any $i=0, \ldots, n$, we have an isomorphism
\[\cH^i(\gL^\cdot(\cF), d_\cdot) \,\, \tilde{\ra} \,\,
i_*(\gL^i(\cE))\] 
of $G$-modules on $\tilde{X}$. 

{\bf Proof}. We subsequently prove stalk-wise that the adjunction
homomorphism 
\[\ker(d_i) \subseteq \gL^i(\cF) \ra i_*i^*(\gL^i(\cF))\]
induces an isomorphism
\[\cH^i(\gL^\cdot(\cF), d_\cdot) \,\, \tilde{\ra}\,\, 
i_*(\gL^i(\cE)).\]
This isomorphism is then compatible with the $G$-actions since this
already holds for the adjunction homomorphism. \\
So let $x \in \tilde{X}$. If $x$ is not contained in $i(X)$, then
$\cI_x = \cO_{\tilde{X},x}$ and the Koszul complex $(\gL^\cdot(\cF)_x,
d_\cdot)$ is exact by Proposition 2.1(c) on p.\ 71 in \cite{FL}.
Thus we have $\cH^i(\gL^\cdot(\cF), d_\cdot)_x = 0 =
i_*(\gL^i(\cE))_x$, as was to be shown.\\
Now let $x\in i(X)$. Then the ideal $I:= \cI_x$ of the noetherian 
local ring $A:= \cO_{\tilde{X},x}$ is generated by a regular sequence,
say of length $d$. Let $H$ be a free $A$-module of rank $d$ and 
$\delta: H \,\,\rightepi\,\, I$ a surjective $A$-homomorphism. Since $F:=
\cF_x$ and $H$ are free, there are homomorphisms $\alpha: H \ra F$ and
$\beta: F\ra H$ such that the following diagrams commute:
\[\begin{array}{ccccc}
H && \stackrel{\alpha}{\ra} && F \\
\\ 
& {\ss \delta} \searrow \phantom{\ss \delta} && 
\phantom{\ss \varepsilon} \swarrow {\ss \varepsilon}\\
\\
&&I
\end{array} 
\qquad {\rm and} \qquad
\begin{array}{ccccc}
F && \stackrel{\beta}{\ra} && H\\
\\
& {\ss \varepsilon} \searrow \phantom{\ss \varepsilon} && 
\phantom{\ss \delta} \swarrow {\ss \delta}\\
\\
&& I.
\end{array}\]
Since $\bar{\delta}: H/IH \ra I/I^2$ is bijective, the composition
$\bar{\beta} \circ \bar{\alpha} : H/IH \ra H/IH$ is bijective, too. 
Using the Nakayama-Lemma, we obtain from this that $\beta \circ 
\alpha: H\ra H $ is bijective. We thus have the decomposition $F =
H \oplus K$, where $H$ is identified with $\alpha(H)$ and where 
$K:= \ker(\beta)$. Since $\varepsilon$ vanishes on $K$, the Koszul
complex $(\gL^\cdot(F), d_\cdot)$ is the tensor product of the
Koszul complex $(\gL^\cdot(H), e_\cdot)$ associated with the 
homomorphism $\delta: H \ra I \subseteq A$ and the complex
$(\gL^\cdot(K), 0)$ with the trivial differential $0$. Thus, we have
by Proposition 2.1(a) on p.\ 71 in \cite{FL}:
\begin{eqnarray*}
\lefteqn{ \cH^i(\gL^\cdot(\cF), d_\cdot)_x \cong H^i(\gL^\cdot(F),
d_\cdot) }\\
&& \cong A/I \otimes \gL^i(K) \cong \gL^i(K/IK) \cong 
\gL^i(\cE_x) \cong i_*(\gL^i(\cE))_x,
\end{eqnarray*}
as was to be shown.

\subsection*{3. $\gS_l$-invariant sections of $\cL^{\boxtimes l}$}

Let $f: X \ra Y$ be a morphism between noetherian schemes, $\cL$ an
invertible $\cO_X$-module, and $l\in \NN$.

{\bf Lemma 3.4}. If $\cL$ is generated by global sections, then 
the $l$-th external tensor power $\cL^{\boxtimes l}$ on $X^l$ is generated
by $\gS_l$-invariant global sections.

{\bf Proof}. Let $p_1, \ldots, p_l: X^l \ra X$ denote the
projections, and let $x \in X^l$. We inductively define sections $s_1, 
\ldots, s_N \in \gG(X,\cL)$ and subsets $M_1, \ldots, M_N$ of 
$I_l = \{1, \ldots, l\}$ as follows: Let $j\in \NN_0$ and let
$s_1, \ldots, s_j$ and $M_1, \ldots, M_j$ be already defined. If 
$M_1 \cup \ldots \cup M_j = I_l$, then we set $N:= j$. Otherwise, 
there exists an $s_{j+1} \in \gG(X,\cL)$ with 
\[M_{j+1} :=
\{i \in I_l \backslash (M_1 \cup \ldots \cup M_j) : 
s_{j+1}(p_i(x)) \not= 0\} \not= \emptyset\] 
by assumption. Now we set
\[s:= \sum_{\gs \in \gS_l/\gS(M_1, \ldots, M_N)} \otimesm_{i\in I_l}
s_{r_{\gs^{-1}(i)}} \in \gG(X^l, \cL^{\boxtimes l})\]
where $r_i := j$ for $i \in M_j$ and where $\gS(M_1, \ldots, M_N)$
denotes the subgroup of those permutations $\gs \in \gS_l$ which
satisfy $\gs(M_j) \subseteq M_j$ for all $j = 1, \ldots, N$. Note
that the sum does not depend on the chosen system of representatives
in $\gS_l$ for $\gS_l/\gS(M_1, \ldots, M_N)$. Thus, $s$ is a
$\gS_l$-invariant section of $\cL^{\boxtimes l}$. Furthermore,
we have 
\[s(x) = \otimesm_{i \in I_l} s_{r_i}(p_i(x)) \not= 0\]
by construction. This proves Lemma 3.4.

\subsection*{4. The $\gS_l$-structure of the conormal sheaf associated
with the diagonal $\gD: X \ra X^l$}

Let $f: X \ra Y$ be a separated smooth morphism between noetherian 
schemes. Let $\cI := \ker(\cO_{X^l} \ra \gD_*(\cO_X))$ denote the
$\gS_l$-stable ideal in $\cO_{X^l}$ associated with the diagonal
$\gD: X \ra X^l$. 

{\bf Lemma 3.5}. We have an isomorphism
\[\cI/ \cI^2 \cong \gO_{X/Y} \otimes \HXl\]
of $\gS_l$-modules on $X$; here, the sheaf $\gO_{X/Y}$ of relative
differentials is considered as a $\gS_l$-module with trivial 
$\gS_l$-action.

{\bf Proof}. The $\cO_X$-module $\HXl$ is free with basis $[i] - 
[i+1]$, $i=1, \ldots, l-1$, (see Notations). For any $i=1, \ldots, 
l-1$, the projection $p_{i,i+1} : X^l \ra X \times X$ defines an
$\cO_X$-homomorphism $p_{i, i+1}^*: \gO_{X/Y} \ra \cI/\cI^2$. We 
define an $\cO_X$-homomorphism 
\[\alpha: \gO_{X/Y} \otimes \HXl \ra \cI/\cI^2\]
by $\alpha(\omega \otimes ([i]-[i+1])) := p_{i,i+1}^*(\omega)$ for
$\omega \in \gO_{X/Y}$ and $i \in \{1, \ldots, l-1\}$. The short 
exact sequences of conormal sheaves associated with the composition
$X \,\, \stackrel{\gD_2}{\longrightarrow} \,\, X\times X \,\,
\stackrel{1\times \gD_2}{\longrightarrow} \,\, X\times X \times X \,\,
\stackrel{1 \times 1\times \gD_2}{\longrightarrow} \,\, \ldots 
\longrightarrow X^l$ (see Proposition 3.4 on p.\ 79 in 
\cite{FL}) show that $\alpha$ is an $\cO_X$-isomorphism. The 
following computation shows that $\alpha$ is compatible with the
$\gS_l$-structures: Let $i \in \{1, \ldots, l-1\}$ and $\gs \in 
\gS_l$. We suppose that $\gs(i+1) > \gs(i)$. (If $\gs(i+1) < \gs(i)$,
one has a similar computation.) Let $p_1, \ldots, p_l: X^l \ra X$
denote the projections. For any local section $x$ in $\cO_X$, we then
have:
\begin{eqnarray*}
\lefteqn{\alpha \gs \left(dx \otimes ([i]-[i+1])\right) = 
\alpha\left(dx \otimes ([\gs(i)]- [\gs(i+1)])\right) }\\
&&= \alpha \left( \sum_{j=\gs(i)}^{\gs(i+1)-1} dx \otimes 
([j]-[j+1])\right)\\
&& = \sum_{j=\gs(i)}^{\gs(i+1)-1}(p_j^*(x) - p_{j+1}^*(x)) +
\cI^2\\
&& = (p_{\gs(i)}^*(x) - p_{\gs(i+1)}^*(x)) + \cI^2 \\
&&= \gs(p_i^*(x) - p_{i+1}^*(x)) + \cI^2\\
&& = \gs \alpha\left(dx \otimes ([i] - [i+1])\right).
\end{eqnarray*}
Since $\gO_{X/Y}$ is locally generated by the total differentials
$dx$, $x\in \cO_X$, this computation shows the required equality
$\gs \alpha = \alpha \gs$. 

\subsection*{5. On the equivariant structure of the ideal $\cI$ in 
$\cO_{X^l}$ associated with the diagonal $\gD: X \ra X^l$}

Let $f: X \ra Y$ be a quasi-projective morphism between noetherian 
schemes (in the sense of \cite{H}, p.\ 103) and $\cL$ an
invertible $\cO_X$-module which is very ample relative to $Y$ (in
the sense of \cite{H}, p.\ 120). Let $\cI:= \ker(\cO_{X^l} \ra
\gD_*(\cO_X))$ denote the $\gS_l$-stable ideal in $\cO_{X^l}$
associated with the diagonal $\gD: X \ra X^l$.

{\bf Proposition 3.6}. There is an epimorphism
\[\alpha: \oplusm^n \cH_{X^l,l} \rightepi \cI \cL^{\boxtimes l}\]
of $C_l$-modules on $X^l$ for some $n > 0$. 

{\bf Proof}. Let $x_0, \ldots, x_r$ be the global sections of $\cL$
associated with an embedding $i: X \hookrightarrow \PP^r_Y$ with
$i^*(\cO(1)) \cong \cL$. For any $j=(j_1, \ldots, j_l) \in 
\{0, \ldots, r\}^l$, let the homomorphism
\[\alpha_j: \cO_{X^l}[I_l] \ra \cL^{\boxtimes l}\]
of $C_l$-modules on $X^l$ be defined by $[i] \mapsto 
c^i(x_{j_1} \otimes \ldots \otimes x_{j_l})$ for $i =1, \ldots, l$. 
By restricting, the map $\alpha_j$ obviously induces a homomorphism
\[\alpha_j: \cH_{X^l,l} \ra \cI \cL^{\boxtimes l}\]
of $C_l$-modules on $X^l$. We set $n:= (r+1)^l = 
\#\{0, \ldots, r\}^l$ and 
\[\alpha:=(\alpha_j)_j : \oplusm^n \cH_{X^l,l} \ra 
\cI \cL^{\boxtimes l}.\]
We are now going to show that $\alpha$ is surjective. For this, we
may assume that $Y =\Spec(A)$ is affine and that $i: X \ra \PP^r_Y$
is a closed immersion. By Corollary 5.16(a) on p.\ 119 in \cite{H}, 
we may furthermore assume that $X= \Proj(S)$ where 
$S= \oplus_{d\ge 0} S_d$ 
is a graded algebra with $S_0 = A$ which is generated by the
elements $x_0, \ldots, x_r \in S_1$. By Exercise 5.11 on p.\ 125 in
\cite{H}, we then have $X^l \cong \Proj(S^l)$ where 
\[S^l := \oplusm_{d\ge 0}  S_d \otimes \ldots \otimes S_d,\] 
and $\cI$ is the
sheaf associated with the homogeneous ideal $I:= \ker(S^l \ra S)$. 
For any $j\in \{0, \ldots, r\}^l$, the map $\alpha_j: 
\cH_{X^l,l} \ra \cI\cL^{\boxtimes l}$ corresponds to the homomorphism
\[\begin{array}{cccl}
\alpha_j:& \oplusm_{i=1}^{l-1} S^l([i] -[i+1]) & \ra & I[1]\\
&[i]-[i+1] &\mapsto & c^i(x_{j_1} \otimes \ldots \otimes
x_{j_l}) - c^{i+1}(x_{j_1} \otimes \ldots \otimes x_{j_l}).
\end{array}\]
Since the open affine subsets $D(x_{k_1} \otimes \ldots \otimes x_{k_l})$, 
$k=(k_1, \ldots, k_l) \in \{0, \ldots, r\}^l$, form a covering of
$\Proj(S^l)$, it suffices to show that, for any $k \in \{0, \ldots,r
\}^l$ and for any homogeneous $f\in I$, the element 
$(x_{k_1} \otimes \ldots \otimes x_{k_l})^N \cdot f$ is contained in
the image of $\alpha = (\alpha_j)_j$ for sufficiently large $N$ (see
also Exercise 5.10 on p.\ 125 in \cite{H}). \\
Now let $(k_1, \ldots, k_r) \in \{0, \ldots, r\}^l$ be fixed. For any
homogeneous $f_1, \ldots, f_l \in S$ of degree $d$, we have:
\begin{eqnarray*} 
\lefteqn{(x_{k_1} \otimes \ldots \otimes x_{k_l})^d \cdot 
(f_1 \otimes \ldots \otimes f_l)}\\
&& = (x_{k_1}^d \otimes \ldots \otimes x_{k_{l-1}}^d \otimes f_l) \cdot
(f_1 \otimes \ldots \otimes f_{l-1} \otimes x_{k_l}^d -
x_{k_l}^d \otimes f_1 \otimes \ldots \otimes f_{l-1})\\
&&\phantom{==} + (x_{k_1} x_{k_l})^d 
\otimes (x_{k_2}^d f_1) \otimes \ldots \otimes 
(x_{k_{l-1}}^d f_{l-2}) \otimes f_{l-1}f_l.
\end{eqnarray*}
We can make a similar computation for
the error term $(x_{k_1}x_{k_l})^d \otimes 
(x_{k_2}^d f_1) \otimes \ldots \otimes (x_{k_{l-1}}^d f_{l-2}) 
\otimes f_{l-1}f_l$ in place of $f_1 \otimes \ldots \otimes f_l$. 
After $(l-1)$ such computations, one sees that, for some sufficiently
large $N$, there are elements $h_1, \ldots, h_l \in S$ which
do not depend on $f_1, \ldots, f_l$ such that the element
\[(x_{k_1} \otimes \ldots \otimes x_{k_l})^N (f_1 \otimes 
\ldots \otimes f_l) - h_1 \otimes \ldots \otimes h_{l-1} \otimes
h_l f_1 \cdot \ldots \cdot f_l\]
is a linear combination of elements of the form $h-c(h)$, $h \in S^l$.
Hence, for any homogeneous $f$ which is contained in the ideal $I$,
the element $(x_{k_1} \otimes \ldots \otimes x_{k_l})^N \cdot f$
is a linear combination of elements of the form $(1-c)(u_1 \otimes 
\ldots \otimes u_l)$, where $u_1, \ldots, u_l$ are monomials in 
the $x_0, \ldots, x_r$ of the same degree. Using the identity 
\[(1-c)(uv) = u((1-c)v) + ((1-c)u)v\] 
(for homogeneous $u,v \in S^l$),
we finally obtain that $(x_{k_1} \otimes \ldots \otimes x_{k_l})^N
\cdot f$ is contained in the image of $\alpha$, as was to be shown.

{\bf Remark 3.7}. If $X= \Spec(B)$ and $Y = \Spec(A)$ are affine, 
then there even exists a homomorphism
\[ \alpha: \oplusm^n \cH_{X^l,l} \ra \cI\]
of $\gS_l$-modules on $X^l$ for some $n > 0$. 

{\bf Proof}. Let $x_1, \ldots, x_n$ be generators of the $A$-algebra
$B$. For any $j=1, \ldots, n$, let the $B^{\otimes l}$-homomorphism
\[\alpha_j: B^{\otimes l}[I_l] \ra B^{\otimes l}\]
be defined by $[i] \mapsto 1 \otimes \ldots \otimes 1 \otimes 
x_j \otimes 1 \otimes \ldots \otimes 1$, where $x_j$ is put 
at the $i$-th place (for any $i=1, \ldots, l$). Obviously, the
map $\alpha_j$ is compatible with the $\gS_l$-structures. As in 
Proposition 3.6, one easily shows that $\alpha := (\alpha_j)_{j=1,
\ldots, n}$ induces an epimorphism $\oplusm^n\cH_{X^l,l} \ra \cI$.

\subsection*{6. Proof of Theorem 3.1}

The element $\lambda_{-1}([\gO_{X/Y} \otimes \HXl])$ is invertible
in $K_0(\gS_l, X)[l^{-1}]/IK_0(\gS_l, X)[l^{-1}]$ 
by Lemma (4.3) in \cite{KoGRR} and by Proposition 3.2.\\
Now let $\cL$ be a very ample invertible $\cO_X$-module relative to
$f$. Let $\cM$ denote the invertible $\gS_l$-module 
$(\cL^{\boxtimes l})^\vee$ on $X^l$. By Lemma 3.4, there is an 
epimorphism
\[\beta: \oplusm^n \cM\,\, \rightepi \,\,\cO_{X^l}\]
of $\gS_l$-modules on $X^l$ for some $n > 0$. The Koszul complex
\[0 \ra \gL^n(\oplusm^n \cM) \ra \ldots \ra \oplusm^n \cM \ra 
\cO_{X^l} \ra 0\]
associated with this epimorphism is exact by Proposition 2.1(c) on p.\
71 in \cite{FL}. Thus we have $(1-[\cM])^n = \lambda_{-1}([\oplusm^n
\cM]) = 0$, i.e., $1-[\cM]$ is nilpotent in $K_0(\gS_l, X^l)$. By
Proposition 3.2, we have
\[\lambda_{-1}([\cM \otimes \cH_{X^l, l}]) = \gt^l(\cM) = 
l + ([\cM]-1) + \ldots +([\cM]^{l-1} -1);\]
thus, the element $\lambda_{-1}([\cM \otimes \cH_{X^l,l}])$ is
invertible in 
$K_0(\gS_l, X^l)[l^{-1}]/IK_0(\gS_l,X^l)[l^{-1}]$.
We set $\cI := \ker(\cO_{X^l} \ra \gD_*(\cO_X))$. By Proposition 3.6, 
there is an epimorphism
\[\oplusm^n(\cM \otimes \cH_{X^l,l}) \rightepi \cI\]
of $C_l$-modules on $X^l$ for some $n >0$. Let the locally free 
$C_l$-module $\cE$ on $X$ be defined by the short exact sequence
\[0 \ra \cE \ra \gD^*\left(\oplusm^n(\cM \otimes \cH_{X^l,l})\right) \ra 
\gD^*(\cI) = \cI/\cI^2 \ra 0.\]
Then we have in $K_0(C_l, X)$ by Lemma 3.5:
\begin{eqnarray*} 
\lefteqn{\lambda_{-1}([\cE]) \cdot \lambda_{-1}([\gO_{X/Y} \otimes
\HXl]) = \lambda_{-1}([\cE]) \cdot \lambda_{-1}([\cI/\cI^2])}\\
&&= \lambda_{-1}\left(\left[\gD^*\left(\oplusm^n(\cM 
\otimes \cH_{X^l,l})\right)\right]\right) = 
\gD^*\left(\lambda_{-1}([\cM \otimes \cH_{X^l,l}])\right)^n.
\end{eqnarray*}
Hence, we have in $K_0(C_l,X^l)[l^{-1}]/IK_0(C_l, X^l)[l^{-1}]$:
\begin{eqnarray*}
\lefteqn{\gD_*\left(\lambda_{-1}([\gO_{X/Y} \otimes
\HXl])^{-1}\right)} \\
&& =\gD_*\left(
\frac{\lambda_{-1}([\cE])}{\gD^*\left(\lambda_{-1}([\cM \otimes
\cH_{X^l,l}])\right)^n } \right)
= \frac{\gD_*(\lambda_{-1}([\cE]))}{\lambda_{-1}\left(\left[\oplusm^n 
(\cM \otimes \cH_{X^l,l})\right]\right)} = 1.
\end{eqnarray*}
Here, the second equality follows from the projection formula; the 
last equality follows from Lemma 3.3 since $\gD$ is regular and all
homology modules, kernels, and images in the Koszul complex
associated with 
$\oplusm^n(\cM \otimes \cH_{X^l,l}) \,\,\rightepi\,\, \cI
\subseteq \cO_{X^l}$
are contained in the category $\cP_\infty(C_l,X^l)$ (defined as in
section 2). 

{\bf Remark 3.8}. If $X$ and $Y$ are affine, then the equality
\[\gD_*\left(\lambda_{-1}([\gO_{X/Y} \otimes \HXl])^{-1}\right) = 1\]
holds even in $K_0(\gS_l, X^l)[l^{-1}]/IK_0(\gS_l, X^l)[l^{-1}]$.

{\bf Proof}. This can be proved as Theorem 3.1 by using 
Remark 3.7 in place of Proposition 3.6.

{\bf Remark 3.9}. If $l$ is a prime, then the ideal $IK_0(C_l,Y)$
in $K_0(C_l,Y)$ is generated by the element $[\cO_Y[I_l]]$.

{\bf Proof}. This follows from the fact that the exterior powers 
$\gL^i(\cO_Y[I_l])$, $i=1, \ldots, l-1$, are free $\cO_Y[C_l]$-modules
(see Proposition 1.1 in \cite{KoCl}).

\bigskip

\section*{\S 4 Riemann-Roch Formulas for Tensor Power Operations}

Let $l \in \NN$, and let $f: X \ra Y$ be a projective morphism as in
section 2. Whilst we considered the $l$-th {\em external} tensor power 
operation for $X$ in section 2, we now consider the $l$-th tensor power 
operation $\tau^l: K_q(X) \ra K_q(\gS_l,X)$. 
The aim of this section is to prove the Riemann-Roch formula
\[\tau^l(f_*(x)) = f_*\left(\lambda_{-1}(T_f^\vee \cdot \HXl)^{-1} \cdot
\tau^l(x)\right)\]
for $x \in K_q(X)$ ($T_f^\vee \in K_0(X)$ denotes the cotangential
element). In general, this formula is valid only in $K_q(C_l,
Y)[l^{-1}]$ modulo a certain subgroup (see Theorem 4.2). It is then 
essentially equivalent to the usual Adams-Riemann-Roch formula for 
$f$ (by Proposition 1.13, Proposition 3.2 and Lemma 4.3). If $f$ is
smooth, we prove Theorem 4.2 by generalizing an idea of Nori (see
\cite{Ra}), and we in particular obtain a new and simple proof of
the usual Adams-Riemann-Roch formula for $f$ in the case $X$ and $Y$ are
$\CC$-schemes. Furthermore, we prove the following strengthened 
versions: If $f$ is a regular closed immersion, the formula above
holds even in $K_q(\gS_l, Y)$ (see Theorem 4.1); here, the multiplier
$\lambda_{-1}(T_f^\vee \cdot \HXl)^{-1}$ equals $\lambda_{-1}
([\cC \otimes \HXl])$ where $\cC$ denotes the conormal sheaf. If $l$
is a prime and $f$ a principal $G$-bundle for some finite group $G$
with $l \notteilt \ord(G)$, then the formula above holds even in 
$K_q(C_l,Y)/[\cO_Y[C_l]] \cdot K_q(C_l,Y)$, i.e., without inverting 
$l$ (see Theorem 4.9); here, the multiplier $\lambda_{-1}(T_f^\vee
\cdot \HXl)^{-1}$ is $1$. In case of $\CC$-schemes, we deduce from
(a $G$-equivariant version of) this result 
a version without denominators of the ($G$-equivariant) 
Adams-Riemann-Roch formula $\psi^l f_* = f_* \psi^l$ (see 
Corollary 4.10).

{\bf Theorem 4.1}. Let $\tilde{X}$ be a noetherian scheme on which
each coherent $\cO_{\tilde{X}}$-module is a quotient of a locally
free $\cO_{\tilde{X}}$-module. Furthermore, let $i:X \hookrightarrow
\tilde{X}$ be a regular closed immersion with the conormal sheaf 
$\cC$. Then the following diagram commutes:
\[\begin{array}{ccc}
K_0(X) & \stackrel{\lambda_{-1}([\cC \otimes \HXl]) \cdot \tau^l}
{\longrightarrow} & K_0(\gS_l, X)\\
\\
{\ss i_*} \downarrow \phantom{\ss i_*} && 
\phantom{\ss i_*} \downarrow {\ss i_*}\\
\\
K_0(\tilde{X}) & \stackrel{\tau^l}{\longrightarrow} & K_0(\gS_l, 
\tilde{X}).
\end{array}\]
If Conjecture' on p.\ 289 in \cite{KoSh} is true, then the
corresponding diagrams of higher $K$-groups commute, too.

{\bf Proof}. First, let $\tilde{X} = \PP_X(\cE \oplus \cO_X)$ where
$\cE$ is a locally free $\cO_X$-module, and let $i:X 
\hookrightarrow \tilde{X}$ be the zero section embedding. Let $\cD$
denote the universal hyperplane sheaf on $\tilde{X}$. For any
$y \in K_0(\tilde{X})$, we then have in $K_0(\gS_l, \tilde{X})$:
\begin{eqnarray*}
\lefteqn{\tau^l(i_*(i^*(y)) = \tau^l(i_*(1) \cdot y) \qquad
\mbox{(projection formula)}}\\
&&= \tau^l(i_*(1)) \cdot \tau^l(y) \qquad \mbox{(Proposition
1.7(a))}\\
&&= \tau^l(\lambda_{-1}([\cD])) \cdot \tau^l(y) \qquad
\mbox{(Lemma 6.2 on p.\ 142 in \cite{FL})}\\
&& = \lambda_{-1}([\cD]) \cdot \lambda_{-1}([\cD \otimes
\cH_{\tilde{X},l}]) \cdot \tau^l(y) \qquad 
\mbox{(Corollary 1.11)}\\
&& = i_*(1) \cdot \lambda_{-1}([\cD \otimes \cH_{\tilde{X}, l}]) \cdot
\tau^l(y) \qquad \mbox{(Lemma 6.2 on p.\ 142 in \cite{FL})}\\
&& = i_*\left(\lambda_{-1}([i^*(\cD \otimes \cH_{\tilde{X},l})]) \cdot
\tau^l(i^*(y))\right) \qquad \mbox{(projection formula)}\\
&& = i_*\left(\lambda_{-1}([\cC \otimes \HXl]) \cdot
\tau^l(i^*(y))\right) \quad
\mbox{(Proposition 3.2(b) on p.\ 78 in \cite{FL})}.
\end{eqnarray*}
Since $i^*: K_0(\PP_X(\cE \oplus \cO_X)) \ra K_0(X)$ is surjective, 
this computation proves Theorem 4.1 for $K_0$ in the case of a 
zero section embedding. The general case follows from this by using
the deformation to the normal cone as, for example, in Theorem 6.3 
on p.\ 142 in \cite{FL}. For higher $K$-groups, Theorem 4.1 can
be proved similarly by using Proposition 1.7(b) in place of 
Proposition 1.7(a). 

Now, let $Y$ be a noetherian scheme on which each coherent
$\cO_Y$-module is a quotient of a locally free $\cO_Y$-module, and
let $f:X \ra Y$ be a projective local complete intersection morphism.
Let $T_f^\vee \in K_0(X)$ denote the cotangential element: For any
decomposition $X \, \, \stackrel{i}{\ra} \,\, \tilde{X} \,\,
\stackrel{p}{\ra} \,\, Y$ of $f$ into a regular closed immersion $i$
and a smooth morphism $p$, we have
\[T_f^\vee = [i^*(\gO_{\tilde{X}/Y})] - [ \cC_{X/\tilde{X}}] \quad {\rm in}
\quad K_0(X)\]
by Proposition 7.1 on p.\ 145 in \cite{FL}. Let $I$ be the ideal
of $K_0(C_l, Y)$ generated by the elements $[\gL^i(\cO_Y[I_l])]$,
$i=1, \ldots, l-1$. By Theorem 3.1, the element $\lambda_{-1}
([\gO_{\tilde{X}/Y} \otimes \cH_{\tilde{X},l}])$ is invertible
in $K_0(C_l, \tilde{X})[l^{-1}]/IK_0(C_l,\tilde{X})[l^{-1}]$. We set
\[\lambda_{-1}(T_f^\vee \cdot [\HXl])^{-1} := 
i^*\left(\lambda_{-1}([\gO_{\tilde{X}/Y} \otimes 
\cH_{\tilde{X},l}])^{-1}\right)
\cdot \lambda_{-1}([\cC_{X/\tilde{X}} \otimes \HXl]) \] 
in $K_0(C_l,X)[l^{-1}]/IK_0(C_l, X)[l^{-1}]$.
This definition does obviously not depend on the chosen decomposition
of $f$ as above.

{\bf Theorem 4.2}. The following diagram commutes:
\[\begin{array}{ccc}
K_0(X)& \stackrel{\lambda_{-1}(T_f^\vee \cdot [\HXl])^{-1} \cdot 
\tau^l}{\longrightarrow} & K_0(C_l,X)[l^{-1}]/IK_0(C_l,X)[l^{-1}]\\
\\
{\ss f_*} \downarrow \phantom{\ss f_*} && 
\phantom{\ss bar{f}_*} \downarrow {\ss \bar{f}_*}\\
\\
K_0(Y) & \stackrel{\tau^l}{\longrightarrow} & 
K_0(C_l,Y)[l^{-1}]/IK_0(C_l,Y)[l^{-1}].
\end{array}\]
The corresponding diagrams of higher $K$-groups commute if one of the
following conditions holds:\\
(a) $f$ is smooth.\\
(b) Conjecture' on p.\ 289 in \cite{KoSh} is true.\\
(c) The number $l$ is a prime, and the 
relation between $\tau^l$ and $\psi^l$ mentioned in 
Proposition 1.13 holds for the higher $K$-groups of $X$ and $Y$. 

{\bf Proof}. If $l$ is  a prime, the assertion for $K_0$ and the
assertion for the higher $K$-groups in case of the condition (c)
follow from the Adams-Riemann-Roch theorem (see Theorem 6.3 on p.\ 142
in \cite{FL} and Theorem (4.5), Remark (4.6)(b) in \cite{KoGRR}) by
using Propositions 1.13 and 3.2.\\
Now, let $f$ be smooth. In this case, we prove Theorem 4.2 for $K_0$
(once more) and for the higher $K$-groups by generalizing Nori's idea (see
\cite{Ra}) to our situation. For this, let $\gD: X \ra X^l$ denote the
diagonal and $f^l: X^l \ra Y$ the projection. By $\tau^l$ we denote
both the tensor power operation $K_q(X) \ra K_q(C_l, X)$ and the
external tensor power operation $K_q(X) \ra K_q(C_l,X^l)$. For any
$x \in K_q(X)$, we then have in $K_q(C_l,X^l)[l^{-1}]/
IK_q(C_l,X^l)[l^{-1}]$:
\begin{eqnarray*}
\lefteqn{\gD_*\left(\lambda_{-1}([\gO_{X/Y} \otimes \HXl])^{-1} \cdot 
\tau^l(x)\right)}\\
&&= \gD_*\left(\lambda_{-1}([\gO_{X/Y} \otimes \HXl])^{-1} \cdot
\gD^*(\tau^l(x))\right) \qquad \mbox{(since $\gD^*(\tau^l(x)) = \tau^l(x)$)}\\
&&= \gD_*\left(\lambda_{-1}([\gO_{X/Y} \otimes \HXl])^{-1}\right) \cdot 
\tau^l(x) \qquad \qquad \mbox{(projection formula)}\\
&&= \tau^l(x) \qquad \qquad \mbox{(Theorem 3.1)}.
\end{eqnarray*}
For any $x\in K_q(X)$, we thus have in $K_q(C_l,Y)[l^{-1}]/
IK_q(C_l,Y)[l^{-1}]$:
\begin{eqnarray*}
\lefteqn{f_*\left(\lambda_{-1} ([\gO_{X/Y} \otimes \HXl])^{-1} \cdot 
\tau^l(x)\right) }\\
&&= f^l_*\gD_*\left(\lambda_{-1}([\gO_{X/Y} \otimes \HXl])^{-1} \cdot
\tau^l(x)\right) \qquad \mbox{(since $f=f^l \circ \gD$)}\\
&& = f^l_*(\tau^l(x)) \qquad \qquad \mbox{(see above)}\\
&& = \tau^l(f_*(x)) \qquad \qquad \mbox{(Theorem 2.1)}.
\end{eqnarray*}
By putting this together with Theorem 4.1, the assertion of Theorem 4.2
for higher $K$-groups also follows in the case (b).

If the canonical map $K_q(Y) \ra K_q(C_l,Y)/IK_q(C_l,Y)$ is 
injective (and if the diagram in Proposition 1.13 commutes for the
higher $K$-groups of $X$ and $Y$), then, conversely, the Adams-Riemann-Roch
theorem obviously follows from Theorem 4.2. If, in addition, $f$ is
smooth, Nori's idea in particular yields a quick proof of the
Adams-Riemann-Roch theorem which, in contrast to \cite{FL}, 
does not use the deformation to
the normal cone and the explicit computations for so-called elementary
embeddings and projections.
The following considerations deal with the question under which
conditions the map $K_q(Y) \ra K_q(C_l,Y)/IK_q(C_l,Y)$ is injective.

{\bf Lemma 4.3}. Let $l$ be  a prime, $q \ge 0$, and $Y$ a scheme
on which $l$ is invertible. Then we have: If the map $K_q(Y) \ra
K_q(C_l,Y)$, $y \mapsto [\HXl]\cdot y$, is injective, then the
canonical map $K_q(Y) \ra K_q(C_l,Y)/IK_q(C_l,Y)$ is injective, too.

{\bf Proof}. By Remark 3.9, the ideal $I$ is generated by the 
element $[\cO[C_l]]$. Furthermore, we have $[\cO_Y[C_l]] \cdot
K_q(C_l,Y) = [\cO_Y[C_l]] \cdot K_q(Y)$ by Frobenius reciprocity.
Finally, we have $K_q(Y) \cap ([\HXl \cdot K_q(Y)) =0$ since
the functor which maps a locally free $C_l$-module $\cF$ to the 
module $\cF^{C_l}$ of $C_l$-fixed elements induces a homomorphism
$K_q(C_l,Y) \ra K_q(Y)$ which is the identity on $K_q(Y)$ and which
vanishes on $[\HXl]\cdot K_q(Y)$. This immediately implies Lemma 4.3.

{\bf Corollary 4.4}. The canonical map $K_q(Y) \ra K_q(C_l,Y)/
IK_q(C_l,Y)$ is injective if $K_q(Y)$ has no $(l-1)$-torsion or if
there exists a primitive $l$-th root of unity on $Y$.

{\bf Proof}. If $K_q(Y)$ has no $(l-1)$-torsion, the composition
\[\begin{array}{ccccc}
K_q(Y) & \ra & K_q(C_l,Y) &  \stackrel{{\rm can}}{\ra} & K_q(Y) \\
y & \mapsto & [\cH_{Y,l}] \cdot y
\end{array}\]
is injective; hence, also the first map is injective. If there exists
a primitive $l$-th root of unity on $Y$, we have $K_q(C_l,Y) \cong
\oplusm^l K_q(Y)$, and each of the non-trivial projections 
$K_q(C_l,Y) \ra K_q(Y)$ is left-inverse to the map $K_q(Y) \ra 
K_q(C_l,Y)$, $y \mapsto [\cH_{Y,l}] \cdot y$. Now, Lemma 4.3 proves
Corollary 4.4.

{\bf Example 4.5}.  Let $l$ be a prime.\\
(a) Let $Y = \Spec(\FF_l)$. Then we have $K_0(C_l,Y)/([\cO[C_l]]) \cong
\ZZ/l\ZZ$. In particular, the canonical map $K_0(Y) \ra 
K_0(C_l,Y)/IK_0(C_l,Y)$ is not injective.\\
(b) Let $Y= \Spec(\CC)$. Then $K_0(C_l,Y)/([\cO_Y[C_l]])$ is 
isomorphic to $\ZZ[\zeta_l]$ where $\zeta_l$ denotes a primitive 
$l$-th root of unity.

Now, let $f:X \ra Y$ be a morphism between noetherian schemes and $G$
a finite group which acts on $X$ by $Y$-automorphisms. We recall the
following definition of \cite{SGA 1}:

{\bf Definition 4.6}. The $Y$-scheme $X$ is called a principal
$G$-bundle, if $f$ is faithfully flat and the morphism
\[G\times X \ra X\times_Y X, \quad (\gamma,x) \mapsto (x,\gamma(x)),\]
is an isomorphism.

Here, we have used the following notation: For any finite set $M$ and
for any scheme $Z$, let $M\times Z$ denote the disjoint union of
$\#M$ copies of $Z$. If $M$ is $G$-set, we will consider $M\times Z$
as a $G$-scheme in the obvious way.

{\bf Remark 4.7}. By Proposition 2.6 on p.\ 115 in \cite{SGA 1}, $X$
is a principal $G$-bundle, if and only if $f$ is finite, $Y=X/G$, and
the inertia groups of all points in $X$ are trivial. In this case, $f$
is \'etale (see Corollaire 2.3 on p.\ 113 in \cite{SGA 1}), and
the exact functor $f_*$ obviously induces 
homomorphisms
\[f_*: K_q(X) \ra K_q(Y) \quad {\rm and} \quad 
f_*: K_q(G,X) \ra K_q(G,Y)\]
for all $q \ge 0$. 

{\bf Example 4.8}. An extension $R \subseteq S$ of noetherian commutative rings
is a Galois extension with Galois group $G$ (in the sense of 
Definition 1.5 on p.\ 2 in \cite{Gre}) if and only if 
$f:\Spec(S) \ra \Spec(R)$ is a principal $G$-bundle. For instance, the
extension of the rings of integers in a Galois extension of number
fields is a Galois extension, if and only if it is unramified (see
Theorem 4.1 on p.\ 18 in \cite{Gre}).

In the following Theorem 4.9, we will consider both the tensor power
operation $\tau^l: K_q(X) \ra K_q(C_l,X)$ and the equivariant 
version $\tau^l: K_q(G,X) \ra K_q(C_l\times G, X)$ which is defined
in the obvious way (see also Remark 1.12). The proof of Theorem 4.9
is similar to the proof of Theorem 4.2; however, we will not use the
results of sections 2 and 3, but more or less only the very definition
of $\tau^l$.

{\bf Theorem 4.9}. Let $l$ be a prime with $l \notteilt \ord(G)$. Let 
$f:X \ra Y$ be a principal $G$-bundle or the canonical projection 
associated with $X= (G/G')\times Y$ where $G'$ is a subgroup of $G$. 
Then the diagram 
\[\begin{array}{ccc}
K_q(H,X) & \stackrel{\tau^l}{\longrightarrow} & 
K_q(C_l\times H,X)/[\cO_X[C_l]] \cdot K_q(C_l \times H,X)\\
\\
{\ss f_*} \downarrow \phantom{\ss f_*} && 
\phantom{\ss \bar{f}_*} \downarrow {\ss \bar{f}_*}\\
\\
K_q(H,Y) & \stackrel{\tau^l}{\longrightarrow} & 
K_q(C_l\times H,Y)/[\cO_Y[C_l]] \cdot K_q(C_l\times H,Y)
\end{array}\]
commutes for all subgroups $H$ of $G$ and all $q \ge 0$. 

{\bf Proof}. The first (big) part of the proof is to 
describe the $(C_l\times G)$-structure of
the fibred product $X^l$ where $C_l$ acts by cyclic permutations
as usual (see Notations) and where $G$ acts diagonally.\\
First, let $f:X \ra Y$ be a principal $G$-bundle. We identify
$G^{l-1}$ with $\{1\} \times G^{l-1} \subseteq G^l$ and define an
action of $C_l$ on $G^{l-1}$ by
\begin{eqnarray*}
\lefteqn{c^i(g_1, \ldots, g_l) :=  g_{i+1}^{-1} \cdot
(g_{i+1}, \ldots, g_l, g_1, \ldots, g_i)}\\
&&=(1,g_{i+1}^{-1}g_{i+2}, \ldots ,g_{i+1}^{-1} g_l, 
g_{i+1}^{-1}, g_{i+1}^{-1}g_2, \ldots,g_{i+1}^{-1} g_i)
\end{eqnarray*}
(for $i \in \{0, \ldots, l-1\}$ and $(g_1, \ldots, g_l) \in G^{l-1}$).
We furthermore define an action of $C_l$ on $G^{l-1} \times X$ by
\[c^i(g,x) := (c^i(g), g_{i+1}(x))\]
(for $g=(g_1, \ldots, g_l) \in G^{l-1}$ and $x\in X$). Finally, we
define a (right) action of $G$ on $G^{l-1} \times X$ by
\[\gamma(g,x) := (\gamma^{-1}g\gamma, \gamma(x))\]
(for $g \in G^{l-1}$, $\gamma \in G$, and $x\in X$). One easily checks
that $G^{l-1}\times X$ together with these actions is a $(C_l\times G)$-scheme
and that the morphism
\[G^{l-1} \times X \ra X^l, \quad ((g_1, \ldots, g_l), x) \mapsto
(g_1(x), \ldots, g_l(x)),\]
is an isomorphism of $(C_l\times G)$-schemes. \\
The action of $G$ on $G^{l-1}$ by conjugation obviously induces an 
action on the set $(G^{l-1}\backslash\{(1, \ldots, 1)\})/C_l$ of 
$C_l$-orbits. We will now show that there is a $G$-stable system of
representatives in $G^{l-1}\backslash\{(1,\ldots,1)\}$ for this
set of orbits. For this, it suffices to show that, for any $g \in
G^{l-1}\backslash \{(1, \ldots,1)\}$, the map 
\[C_l \times \{\gamma^{-1}g\gamma : \gamma \in G\} \ra 
G^{l-1}, \quad (c^i,h) \mapsto c^i(h),\]
is injective. So let $i,j \in \{0, \ldots, l-1\}$ and $\eta, \gt \in
G$ with $c^i(\eta^{-1}g \eta) = c^j(\gt^{-1}g\gt)$. We may assume that
$j=0$ and $\eta = 1$. Then we have $c^i(g) = \gt^{-1} g \gt$, hence
$c^{ir}(g) = \gt^{-r} g \gt$ for all $r\ge 0$; i.e., the
orbit of $g$ under the subgroup of $G$ generated by $\gt$ is the set
$N:= \{c^0(g), c^i(g), c^{2i}(g), \ldots\}$. Thus, 
the number $\#N$ of elements in $N$ divides the
order of $\gt$. Since 
$\#N \in \{1,l\}$ and since $l\notteilt \ord(G)$, we have $\#N =1$,
i.e., $c^i(g) = g$. This implies $i=0$ because, otherwise, $g$
would be fixed under the action of $C_l$ and  would hence be of the
form $(\gamma^0, \gamma^1, \ldots, \gamma^{l-1})$ for some $\gamma \in
G$ with $\gamma^l =1$ which means that $g$ would be equal to 
$(1, \ldots, 1)$ since $l \notteilt \ord(G)$. \\
Now, let $M\subseteq G^{l-1} \backslash \{(1, \ldots, 1)\}$ be a 
$G$-stable system of representatives for 
\[(G^{l-1} \backslash \{(1, \ldots, 1\})/C_l.\] 
We consider $C_l\times M \times X$ as 
a $(C_l \times G)$-scheme by virtue of 
\[(c^i, \gamma)(c^j,g,x):=(c^{i+j}, \gamma^{-1}g \gamma, \gamma(x))\]
(for $i,j \in \{0, \ldots, l-1\}$, $\gamma \in G$. $g\in M$,
and $x\in X$). Then the morphism
\[C_l \times M \times X \ra 
\left(G^{l-1} \backslash \{(1, \ldots,1)\}\right) \times X, \quad
(c^i, g, x) \mapsto (c^i(g), g_i(x)),\]
is obviously an isomorphism of $(C_l\times G)$-schemes. Altogether,
we now have the following isomorphisms of $(C_l\times G)$-schemes:
\[X^l \cong \gD(X) \coprod (G^{l-1} \backslash \{(1, \ldots, 1)\}) 
\times X \cong \gD(X) \coprod C_l \times M \times X\]
where $\gD:X\ra X^l$ denotes the diagonal.\\
Now, let $X=(G/G')\times Y$ where $G'$ is a subgroup of $G$. Let 
$\gD: G'\backslash G \ra (G'\backslash G)^l$ denote the diagonal. 
Similarly as above, we show that, for any $g \in (G'\backslash G)^l
\,\backslash \,\gD(G'\backslash G)$, the map
\[C_l \times gG \ra (G'\backslash G)^l, \quad 
(c^i,h) \mapsto c^i(h),\]
is injective. So let $i,j \in \{0, \ldots, l-1\}$ and 
$\eta, \gt \in G$ with $c^i(g\eta)= c^j(g\gt)$. We may assume that
$j=0$ and $\eta = 1$. Then we have $c^i(g) = g\gt$, hence
$c^{ir}(g)= g\gt^r$ for all $r\ge 0$; i.e., the orbit
of $g$ under the subgroup of $G$ generated by $\gt$ is the set
$N:=\{c^0(g), c^i(g), c^{2i}(g), \ldots\}$. Thus, the 
number $\#N$ of elements in $N$ divides the order of $\gt$. Since
$\# N \in \{1,l\}$ and since $l \notteilt \ord(G)$, we have $\#N =1$,
i.e., $c^i(g) = g$. This implies $i=0$ because, otherwise, we
would have $g \in \gD(G'\backslash G)$. \\
As above, this implies that there exists a $G$-stable subset 
$M\subseteq (G'\backslash G)^l \,\backslash\, \gD(G'\backslash G)$
such that the map
\[C_l \times M \ra (G'\backslash G)^l \,\backslash \,\gD(G'\backslash G),
\quad (c^i, h) \mapsto c^i(h),\]
is bijective. Hence, we have isomorphisms
\[X^l \cong (G'\backslash G)^l \times Y \cong 
\gD(X) \coprod C_l \times M \times Y\]
of $(C_l\times G)$-schemes.\\
Now, let $f$ be a principal $G$-bundle or $X=(G'\backslash G) \times
Y$. If follows for instance from the established decompositions of 
$X^l$ that, for any locally free $\cO_X$-module $\cE$, the 
direct image $\gD_*(\cE)$ is a locally free $\cO_{X^l}$-module. The functor
$\gD_*$ thus induces a homomorphism
\[\gD_*: K_0(C_l \times H, X) \ra K_0(C_l \times H,X^l)\]
for all subgroups $H$ of $G$. Furthermore, it follows from the 
established decompositions of $X^l$ that the homomorphism
$\Ind_{\{1\}}^{C_l}: K_0(H,X^l \backslash \gD(X)) \ra 
K_0(C_l \times H, X^l \backslash \gD(X))$ is surjective. This implies
that the homomorphism
\[\gD_*: K_0(C_l\times H, X)/([\cO_X[C_l]]) \ra 
K_0(C_l \times H, X^l)/\Ind_{\{1\}}^{C_l}(K_0(H,X^l))\]
is an isomorphism with 
\[\gD_*(1) = 1.\]
Now, let $f^l: X^l \ra Y$ denote the projection. By the projection
formula and the base change isomorphism (see Proposition 9.3 
on p.\ 255 in \cite{H}),
we have for all $x\in K_q(H,X)$:
\[\tau^l(f_*(x)) = f^l_*(\tau^l(x)) \quad {\rm in} \quad
K_q(C_l\times H, Y).\]
Furthermore, we obviously have:
\[f^l_*\left(\Ind_{\{1\}}^{C_l}(K_0(H,X^l)) \cdot
K_q(C_l \times H, X^l)\right) \subseteq 
[\cO_Y[C_l]]\cdot K_q(C_l \times H, Y).\]
As in the proof of Theorem 4.2, we now have:
\[\tau^l(f_*(x)) = f^l_*(\tau^l(x)) = f^l_*(\tau^l(x) \cdot
\gD_*(1))
= f^l_*(\gD^*(\tau^l(x))) = f_*(\tau^l(x))\]
in $K_q(C_l\times H,Y)/[\cO_Y[C_l]] \cdot K_q(C_l\times H,Y)$.
This proves Theorem 4.9.

{\bf Corollary 4.10}. Let $Y$ be a $\ZZ[l^{-1}][\zeta_l]$-scheme. 
Then we have for all $x\in K_0(H,X)$:
\[\psi^l(f_*(x)) = f_*(\psi^l(x)) \quad {\rm in} \quad
K_0(H,Y).\]
The corresponding assertion for higher $K$-groups is true, if
the relation between $\tau^l$ and $\psi^l$ mentioned in Proposition
1.13 holds for the (higher equivariant) $K$-groups of $X$ and $Y$. 

{\bf Proof}. This follows from Theorem 4.9 and obvious generalizations
of Proposition 1.13 and Corollary 4.4.

{\bf Remark 4.11}. Let $G'$ be a subgroup of $G$ and $Y$ a 
noetherian scheme. In section 6 in \cite{KoGRR}, we have raised the
following question: Under which conditions does the induction formula
\[\psi^l(\Ind_{G'}^G(x)) = \Ind_{G'}^G(\psi^l(x))\]
hold in $K_q(G,Y)$ for $x\in K_q(G',Y)$? In Theorem (6.2) in 
\cite{KoGRR}, we have proved that, for all $l$, this formula holds in
a certain completion of the group 
$K_q(G,Y)[l^{-1}]$. From Corollary 4.10 it
now follows as in Theorem (6.2) in \cite{KoGRR} that this formula
holds already in $K_q(G,Y)$ if $l$ is a prime with $l \notteilt \ord(G)$
and $Y$ is a $\ZZ[l^{-1}][\zeta_l]$-scheme (and if, in case $q \ge 1$,
the relation between $\psi^l$ and $\tau^l$ mentioned in Proposition
1.13 holds). 

\bigskip

\section*{\S 5 Speculations in Characteristic $p$}

For any scheme $Z$ of characteristic $p$, the $p$-th Adams operation
$\psi^p$ on $K_0(Z)$ equals the pull-back homomorphism $F_Z^*$ 
associated with the absolute Frobenius morphism $F_Z$. This 
well-known fact (e.g., see Proposition 2.15 on p.\ 64 in
\cite{KoAdPro}) should be considered as a substitute for the relation
``$\psi^p = \tau^p$'' established in Proposition 1.13 which, in
characteristic $p$, still holds but is extremely weak (cf.\ 
Example 4.5(a)). In this section, we will formulate and investigate
an assertion (see Question 5.2) which may be considered as a
substitute for Theorem 3.1. In doing so, we will find many
analogies between the previous sections and the considerations in
this section.

So let $p$ be a prime and $f: X\ra Y$ a morphism between noetherian
$\FF_p$-schemes $X$ and $Y$. As a substitute for the external
tensor power operation $\tau^l: K_0(X) \ra K_0(\gS_l, X^l)$,
we consider the pull-back homomorphism $F_Y^*: K_0(X) \ra 
K_0(X_Y)$ where $X_Y$ and $F_Y:X_Y \ra X$ are defined by the 
following cartesian commutative diagram:
\[\begin{array}{ccc}
X_Y & \stackrel{F_Y}{\longrightarrow} & X \\
\\
{\ss f_Y} \downarrow \phantom{\ss f_Y} && 
\phantom{\ss f} \downarrow {\ss f}\\
\\
Y & \stackrel{F_Y}{\longrightarrow} & Y.
\end{array}\]
The following proposition may be viewed as an analogue of Theorem 2.1:

{\bf Proposition 5.1}. Let $Y$ be a noetherian scheme  on which each
coherent $\cO_Y$-module is a quotient of a locally free
$\cO_Y$-module. Furthermore, let $f:X \ra Y$ be a flat projective
local complete intersection morphism. Then the following diagram 
commutes:
\[\begin{array}{ccc}
K_0(X) & \stackrel{F_Y^*}{\longrightarrow} & K_0(X_Y) \\
\\
{\ss f_*} \downarrow \phantom{\ss f_*} && 
\phantom{\ss (f_Y)_*} \downarrow {\ss (f_Y)_*}\\
\\
K_0(Y) & \stackrel{F_Y^*}{\longrightarrow} & K_0(Y).
\end{array}\]

{\bf Proof}. This follows from the excess intersection formula (see
Theorem 1.3 on p.\ 155 in \cite{FL}). Note that the excess conormal sheaf
in our situation is trivial since $f$ is flat.

As a substitute for the diagonal $\gD: X \ra X^l$, we consider the
relative Frobenius morphism $F: X \ra X_Y$ which is defined 
by the following commutative diagram:
\[\begin{array}{ccccc}
X\\
&& {\ss F_X}\\
& \phantom{\ss F} \searrow {\ss F} & \\
\\
&& X_Y & \stackrel{F_Y}{\ra} & X\\
&{\ss f} \\
&& \phantom{\ss f_Y} \downarrow {\ss f_Y} && \phantom{\ss f}
\downarrow {\ss f}\\
\\
&&Y & \stackrel{F_Y}{\ra} & Y.
\end{array}\]

Now, let $f$ be smooth. Then, a theorem of Kunz (e.g., see
Theorems 15.7 and 15.5 in \cite{KuKD}) says that $F$ is
(finite and) flat. This fact may be considered as an analogue of the
fact that the closed immersion $\gD$ is regular if $f$ is smooth. 
As a substitute for the push-forward homomorphism $\gD_*$, we in
particular have the push-forward homomorphism
\[F_*: K_0(X) \ra K_0(X_Y)\]
(induced by the exact functor $F_*$). As in 
section 3, let $\gt^p(\gO)^{-1} \in K_0(X)[p^{-1}]$ denote the
inverse of the Bott element $\gt^p(\gO)$ associated with the module
$\gO := \gO_{X/Y}$ of relative differentials. In view of Proposition 
3.2, the following question
may be considered as an analogue of Theorem 3.1:

{\bf Question 5.2}. Do we have $F_*(\gt^p(\gO)^{-1}) = 1$ in 
$K_0(X_Y)[p^{-1}]$?

The following proposition may be considered as an analogue of the
equality
\[\lambda_{-1}\left(\left[\oplusm^n
(\cM \otimes \cH_{X^l,l})\right]\right) \cdot
\left(\gD_*(\gt^p(\gO)^{-1})-1\right) =0\]
established in the proof of Theorem 3.1.

{\bf Proposition 5.3}. The element $\lambda_{-1}([F_Y^*(\gO)]) = 
\lambda_{-1}([\gO_{X_Y/Y}]) \in K_0(X_Y)$ annihilates the 
difference $F_*(\gt^p(\gO)^{-1}) -1$. 

{\bf Proof}. Using the splitting principle (e.g., see section (2.5) in
\cite{KoARR}), one easily shows that, for any locally free 
$\cO_Y$-module $\cF$, the equality
\[\lambda_{-1}([F_X^*(\cF)]) \cdot \gt^p(\cF)^{-1} = 
\lambda_{-1}([\cF]) \quad {\rm in} \quad K_0(X_Y)[p^{-1}]\]
holds. Thus we have:
\begin{eqnarray*}
\lefteqn{\lambda_{-1}([F_Y^*(\gO)]) \cdot F_*(\gt^p(\gO)^{-1})}\\
&&= F_*\left(\lambda_{-1}([F_X^*(\gO)]) \cdot \gt^p(\gO)^{-1}\right) \\
&& = F_*(\lambda_{-1}([\gO]))\\
&& = \lambda_{-1}([F_Y^*(\gO)]) \quad {\rm in} \quad 
K_0(X_Y)[p^{-1}].
\end{eqnarray*}
Here, the first equality follows from the projection formula and the
last equality follows from the Cartier isomorphism (see Theorem (7.2) 
on p.\ 200 in \cite{Ka}). This proves Proposition 5.3.

In contrast to the element $\lambda_{-1}\left(\left[\oplusm^n(\cM \otimes
\cH_{X^l,l})\right]\right)$, the element $\lambda_{-1}([F_Y^*(\gO)])$ is not
invertible, it is even nilpotent. Thus, Proposition 5.3 does not imply the
equality $F_*(\gt^p(\gO)^{-1}) =1$.

In addition to the previous assumptions, we now suppose that each
coherent $\cO_Y$-module on $Y$ is a quotient of a locally free 
$\cO_Y$-module and that $f$ is projective. We set
\[{\rm Num}_0(f) := \{y \in K_0(X_Y)[p^{-1}]: (f_Y)_*( y \cdot 
F_Y^*(x)) = 0 \mbox{ for all } x \in K_0(X)[p^{-1}]\}.\]
The notation ${\rm Num}_0(f)$ should remind the reader of the group
of cycles numerically equivalent to zero which is defined in a similar
way.

{\bf Example 5.4}. \\
(a) Let $f: \PP(\cE) \ra Y$ be the canonical projection of the 
projective space bundle associated with a locally free $\cO_Y$-module
$\cE$. Then we have ${\rm Num}_0(f)=0$. This follows from the
projective space bundle theorem (see Theorem 2.3 on p.\ 115 in
\cite{FL}).\\
(b) Let $f: X \ra \Spec(\FF_p)$ be a smooth projective curve. Then
we have ${\rm Num}_0(f) \otimes \QQ =0$. This follows from the 
Riemann-Roch theorem for curves (e.g., see Theorem 1.3 on p.\ 295 in
\cite{H}) because $K_0(X) \cong {\rm Pic}^0(X) \oplus \ZZ \oplus \ZZ$.

The following proposition may be considered as an analogue of the
fact that Theorem 3.1 implies the Adams-Riemann-Roch formula for 
$\CC$-schemes (see section 4). This proposition also clarifies 
to what extent, conversely, the Adams-Riemann-Roch formula implies the
equality $F_*(\gt^p(\gO)^{-1})=1$. 

{\bf Proposition 5.5}. The difference $F_*(\gt^p(\gO)^{-1}) -1$ 
is contained in ${\rm Num}_0(f)$ if and only if, for all $x \in
K_0(X)$, the Adams-Riemann-Roch formula
\[\psi^p(f_*(x)) = f_*(\gt^p(\gO)^{-1} \cdot \psi^p(x)) \quad
{\rm in} \quad K_0(Y)[p^{-1}]\]
holds.

{\bf Proof}. For all $x\in K_0(X)$, we have
\[\psi^p(f_*(x)) = F_Y^*(f_*(x)) = (f_Y)_*(F_Y^*(x))\quad {\rm in}
\quad K_0(Y)\]
by Proposition 5.1. Furthermore, using the projection formula, we
obtain:
\begin{eqnarray*}
\lefteqn{f_*(\gt^p(\gO)^{-1} \cdot \psi^p(x)) }\\
&& = (f_Y)_* \left(F_*\left(\gt^p(\gO)^{-1} \cdot 
F^*(F_Y^*(x))\right)\right)\\
&& = (f_Y)_* \left(F_*(\gt^p(\gO)^{-1}) \cdot F_Y^*(x)\right)
\quad {\rm in} \quad K_0(Y)[p^{-1}].
\end{eqnarray*}
This immediately implies Proposition 5.5.

The following remark may be considered as an analogue of the equality
$\gD_*(1) = 1$ established in the proof of Theorem 4.9.

{\bf Remark 5.6}. Let $f: X \ra Y$ be an \'etale morphism. Then 
we have $\gO = 0$, and $F$ is an isomorphism (see Lemma 3.7 on p.\ 163
in \cite {FK}). Hence, we have:
\[F_*(\gt^p(\gO)^{-1}) = F_*(1) = 1 \quad {\rm in} \quad
K_0(X_Y).\]
As in Proposition 5.5, a version without denominators of the
Adams-Rie\-mann-Roch formula $\psi^p f_* = f_* \psi^p$ can be
deduced from this.

In the following example, we describe what the equality
$F_*(\gt^p(\gO)^{-1})=1$ means in case of a smooth
curve.

{\bf Example 5.7}. Let $f: X \ra \Spec(\FF_p)$ be a smooth connected curve. 
Then we have
\[\gt^p(\gO)^{-1} = \frac{2}{p} - \frac{1}{p^2} \gt^p(\gO) \quad
{\rm in} \quad K_0(X)[p^{-1}]\]
since the square of the augmentation ideal $\ker(\rank: K_0(X) \ra
\ZZ)$ vanishes. For instance, we have in the case $p=2$:
\[\gt^2(\gO)^{-1} = \frac{1}{4}(3 - \gO).\]
Now, let $X=\Spec(R)$ be affine. Using the Cartier isomorphism (see 
the proof of Proposition 5.3), one easily shows that, in the case
$p=2$, Question 5.2 is equivalent to the following question: Is 
there an $r\ge 0$ such that 
\[\oplusm^{2^r}\left( F_*(R) \oplus F_*(R)\right) \cong 
\oplusm^{2^r} (R\oplus R \oplus R \oplus \gO)?\]
Furthermore, if $X$ is an elliptic curve, this is equivalent to the 
following question: Is there an $r\ge 0$ such that the direct sum of
$2^r$ copies of $F_*(R)$ is $R$-free?

\bigskip

Department of Mathematics, University of Illinois at Urbana-Champaign,
1409 West Green Street, Urbana, IL 61801, USA; {\em e-mail:}
koeck@math.uiuc.edu (until end of March 1998).

Mathematisches Institut II der Universit\"at Karlsruhe, 
D-76128 Karlsruhe, Germany; {\em e-mail:} 
Bernhard.Koeck@math.uni-karlsruhe.de.

\end{document}